\documentclass[12pt]{article}
\usepackage{array,CJK}
\usepackage{rotating}
\usepackage{lscape}
\usepackage{latexsym,epsfig}
\usepackage{amssymb,amsmath,amsthm}
\usepackage{bm}
\usepackage{cases}
\usepackage[T1]{fontenc}
\usepackage[utf8]{inputenc}
\usepackage{graphicx} 
\usepackage{mwe}
\usepackage{subfig}
\usepackage{authblk}

\newcommand{\al}{\alpha}

\newcommand{\be}{\beta}
\renewcommand{\th}{\theta}
\newcommand{\Th}{\Theta}
\newcommand{\la}{\lambda}

\newcommand{\sig}{\sigma}
\newcommand{\womega}{\widetilde{\omega}}
\newcommand{\rme}{{\rm e}}

\newcommand{\bfN}{{\bf N}}
\newcommand{\bfR}{{\bf R}}

\newcommand{\bfZ}{{\bf Z}}

\newcommand{\bk}{{\bf k}}
\newcommand{\bp}{{\bf p}}

\newcommand{\by}{{\bf y}}

\newcommand{\bq}{{\bf q}}
\newcommand{\bv}{{\bf v}}
\renewcommand{\d}{\displaystyle}

\usepackage{color}

\title{Dispersion relations of periodic quantum graphs associated with Archimedean tilings (II)}
\author[1]{Yu-Chen Luo}
\author[1,2] {Eduardo O.\ Jatulan}
\author[1] {Chun-Kong Law}
 \affil[1]{ \footnotesize{Department of Applied Mathematics, National Sun Yat-sen University,
Kaohsiung, Taiwan 80424. Email: leoredro@gmail.com, law@math.nsysu.edu.tw}}
\affil[2]{ \footnotesize{Institute of Mathematical Sciences and Physics, University of the Philippines Los Ba\~{n}os, Philippines  4031. Email: eojatulan@up.edu.ph}}

\begin{document}
\maketitle
 \begin{abstract}
  We continue the work of a previous paper \cite{LJL} to derive the dispersion relations of the  periodic quantum graphs associated with the remaining 5 of the 11 Archimedean tilings, namely the truncated hexagonal tiling $(3,12^2)$, rhombi-trihexagonal tiling $(3,4,6,4)$, snub square tiling $(3^2,4,3,6)$, snub trihexagonal tiling $(3^4,6)$, and truncated trihexagonal tiling $(4,6,12)$. The computation is done with the help of the symbolic software Mathematica.
  With these explicit dispersion relations, we perform more analysis on the spectra.
\\[0.1in]
 \end{abstract}
 \vskip1in
 \newpage
 \section{Introduction}
 On the plane, there are plenty of ways to form a tessellation with different patterns \cite{GS}. A tessellation is formed
 by regular polygons having the same edge length such that these regular polygons will fit around at
 each vertex and the pattern at every vertex is isomorphic, then it is  called a \textit{uniform tiling}. There are 11 types
 of uniform tilings, called  Archimedean tilings in literature, as shown in Table \ref{fig5.1}. The periodic quantum graphs associated with these
 Archimedean tilings constitute good mathematical models for graphene \cite{ALM,KP} and related allotropes \cite{EI}.  In a previous paper \cite{LJL}, we employ the Floquet-Bloch
 theory, and the characteristic function method for quantum graphs to derive the dispersion relations of periodic quantum graphs
 associated with 4 kinds of Archimedean tilings. As the dispersion relations concerning with square tiling and hexagonal tiling
 have been found earlier \cite{KL,KP}, there are 5 remaining Archimedean tilings: the truncated hexagonal tiling $(3,12^2)$, rhombi-trihexagonal tiling $(3,4,6,4)$, snub square tiling $(3^2,4,3,6)$, snub trihexagonal tiling $(3^4,6)$, and truncated trihexagonal tiling $(4,6,12)$.  The objective of this paper is to derive the dispersion relations associated with
 them, and explore some of the implications.

 Given an infinite graph $G=E(G)\cup V(G)$ generated by an Archimedean tiling with identical edgelength $a$, we let $H$ denote a Schr\"odinger operator on , i.e.,
 $$
 H \, \by(x)=-\frac{d^2}{d x^2}\by(x)+\bq(x)\, \by(x),
 $$
  {where $ \bq\in L^2_{loc}(G)$ is periodic on the tiling (explained below), and} the domain $D(H)$ consists of all  {admissible} functions $\by(x)$ (union of functions $y_e$ for each edge $e\in E(G)$) on $G$ in the sense that
 \begin{enumerate}
 \item[(i)] $\displaystyle   y_e\in {\cal H}^2(e)$ for all $e\in E(G)$;
 \item[(ii)] $\displaystyle   \sum_{ {e\in E(G)}} \| y_e\|^2_{{\cal H}^2(e)} <\infty$;
 \item[(iii)] Neumann vertex conditions (or continuity-Kirchhoff conditions at vertices), i.e., for any vertex $\bv\in V(G)$,
 $$
 y_{e_{i}}(\bv)=y_{e_{j}}(\bv)\qquad \mbox{ and }\qquad \sum_{ {e\in E_\bv(G)}} y_e'(\bv)=0.
 $$
 Here $y_e'$ denotes the directional derivative along the edge $e$ from $\bv$, $E_\bv(G)$ is the set of edges adjacent to $\bv$,  and $\| \cdot\|^2_{{\cal H}^2(e)}$ denotes the Sobolev norm of 2 distribution derivatives.
 \end{enumerate}
 The defined operator $H$ is well known to be self-adjoint. As $H$ is also periodic, we may apply
 the Floquet-Bloch theory \cite{BK,E73,RS} in the study of its spectrum.
  {Let $\vec{k}_1,\vec{k}_2$ be two linear independent vectors in $\bfR^2$. Define $\vec{\bk}=(\vec{k}_1,\vec{k}_2)$, and $\bp=(p_1,p_2)\in \bfZ^2$.
 For any set $S\subset G$, we define $\bp \circ S$ to be an action on $S$ by a shift of $\bp\cdot\bk=p_1 \vec{k}_1+p_2 \vec{k}_2$.
 A compact set $W\subset G$ is said to be a fundamental domain if \
 $$
 G=\bigcup\{ \bp\circ W:\ \bp\in \bfZ^2\},
 $$
 and for any different $\bp,\bp'\in \bfZ^2$, $(\bp\circ W)\cap (\bp'\circ W)$ is a finite set in $G$. The potential function $\bq$ is said to be periodic if
 $\bq(x+\bp\cdot \vec{\bk})=\bq(x)$ for all $\bp\in \bfZ^2$ and all $x\in G$.}
  Take the quasi-momentum $\Theta=(\theta_1,\theta_2)$ in the Brillouin zone
 $B=[-\pi,\pi]^2$. Let $H^{\Th}$ be the Bloch  Hamiltonian that acts on  $L^2(W)$, and the dense domain $D(H^{\Th})$ consists of admissible functions $\by$ which satisfy the Floquet-Bloch condition
 \begin{equation}
 \by(x+\bp\cdot \vec{\bk})= \rme^{i(\bp\cdot\Theta)}\by(x),
 \label{eq1.02}
 \end{equation}
 for all $\bp\in \bfZ^2$ and all $x\in G$.   {Such functions are uniquely determined by their restrictions on the fundamental domain $W$. Hence for fixed $\Th$,}
 the operator $H^\Th$ has purely discrete spectrum $\sig(H^\Th)=\{\la_j(\Th):\ j\in \bfN\}$, where
 $$
  {\la_1(\Th)\leq \la_2(\Th)\leq \cdots \leq \la_j(\Th)\leq\cdots},\quad \mbox{ and } \la_j(\Th)\rightarrow\infty\mbox{ as } j\to\infty.
 $$

  By an analogous argument as in \cite[p.291]{RS}, there is a unitary operator $U:\ L^2(G)\rightarrow L^2(B,\, L^2(W))$ such that
 $$
  UHU^*=\int_B^\oplus H^{\Th}\, d\Th.
 $$
 That is, for any $f\in L^2(B,\, L^2(W))$, $UHU^*\, f(\Th)=H^\Th\, f(\Th)$.
 Hence by \cite{BK,RS},
 $$
  \bigcup \{\sig(H^{\Th}):\ \Th\in [-\pi,\pi]^2\} =\sigma(UHU^*)=\sig(H).
 $$
 Furthermore,  {it is known that  singular continuous spectrum is absent in $\sig(H)$ \cite[Theorem 4.5.9]{K93}}.  These form the basis of our work.

  If one can derive the dispersion relation, which relates the energy levels $\la$ as a function of the quasimomentum  $\Th=(\theta_1,\theta_2)$.
  The spectrum of the operator $H$ will then be given by the set of roots (called Bloch variety or analytic variety) of the dispersion relation.

 On the interval $[0,a]$, we let $C(x,\rho)$ and $S(x,\rho)$ are the solutions of
 $$
 -y''+q y=\lambda y
 $$
 such that $C(0,\rho)=S'(0,\rho)=1$, $C'(0,\rho)=S(0,\rho)=0$.  In particular,
 \begin{eqnarray*}
 C(x,\rho)&=&\displaystyle\cos(\rho x)+\frac1{\rho}\int^x_0\sin(\rho (x-t))q(t)C(t, \rho)dt,\\
 S(x, \rho)&=&\displaystyle \frac{\sin(\rho x)}{\rho}+\frac1{\rho}\int^x_0\sin(\rho (x-t))q(t)S(t, \rho)dt.
 \end{eqnarray*}
 Furthermore it is well known that the Lagrange identity $C S'-S C'=1$, and when $q$ is even, $C(a,\rho)=S'(a,\rho)$ \cite[p.8]{MW}.

  So for the periodic quantum graphs above, suppose the edges $\{ e_1,\ldots,e_I\}$ lie in a typical fundamental domain $W$, while
  $(q_1,\ldots,q_I)$ are the potential functions acting on these edges.
 Assume also that the potential functions $q_i$'s are identical and even.
 Then the dispersion relation can be derived using a characteristic function approach  \cite{LLW,LJL}.
  \newtheorem{th1.0}{Theorem}[section]
 \begin{th1.0}[\cite{KL,KP,LLW,LJL}]
 \label{th1.0}
  Assume that all the $q_j$'s are identical (denoted as $q$), and even. We also let $\theta_1,\ \theta_2\in[-\pi,\pi]$. The dispersion relation of the
  periodic quantum graph associated with each Archimedean tiling is given by the following.
 \begin{enumerate}
 \item[(a)] For $H_S$ associated with square tiling,
 $$
 S(a,\rho)^2\, (S'(a,\rho)^2-\cos^2(\frac{\theta_1}{2})\cos^2(\frac{\theta_2}{2}))=0.
 $$
 \item[(b)] For  $H_H$ associated with hexagonal tiling,
 $$
 S(a,\rho)^2\,\left( 9 S'(a,\rho)^2-1-8\cos(\frac{\theta_1}{2})\cos(\frac{\theta_2}{2})\cos(\frac{\theta_1-\theta_2}{2})\right)=0.
 $$
 \item[(c)] For  $H_{T}$ associated with triangular tiling,
 $$
 S(a,\rho)^2 \left(3S'(a,\rho)+1-4 \cos(\frac{\theta_1}{2})\cos(\frac{\theta_2}{2})\cos(\frac{\theta_2-\theta_1}{2})\right)=0.
 $$
 \item[(d)] For  $H_{eT}$ associated with elongated triangular tiling,
  $$
  S(a,\rho)^3\, \{25(S'(a,\rho))^2-20\cos\theta_1S'(a,\rho)-8\cos(\frac{\theta_1}{2})\cos(\frac{\theta_2}{2})\cos(\frac{\theta_1-\theta_2}{2})+4\cos^2\theta_1-1\}=0.
 $$
 \item[(e)] For  $H_{trS}$ associated with truncated square tiling,
 $$
 S(a,\rho)^2\left\{ 81 (S'(a,\rho))^4-54 (S'(a,\rho))^2-12S'(a,\rho)\, (\cos\theta_1+\cos\theta_2)+1-4\cos\theta_1\, \cos\theta_2\right\} = 0.
 $$
 \item[(f)] For  $H_{TH}$ associated with trihexagonal tiling,
 $$
 S(a,\rho)^3(2S'(a,\rho)+1)\, \left(2 (S'(a,\rho))^2-S'-cos\frac{\theta_1}{2}\cos\frac{\theta_2}{2}\cos\frac{\theta_1-\theta_2}{2})\right)=0.
 $$
 \end{enumerate}
 \end{th1.0}
We note that there was a typo error in \cite[Theorem 1.1(d)]{LJL} about this relation. In fact, in \cite{LLW,LJL}, the dispersion relations were derived \underline{without} the restriction that $q_j$'s are identical and even, although the equations are more complicated.
 Also observe that all the above dispersion relations are of the form $S^i(S')^j\, p(S',\th_1,\th_2)=0$, where $p$ is a polynomial in $S'(a,\rho)$, with
  $p(1,0,0)=0$. Note that the term $S^i(2S'+1)^j$ determines the point spectrum $\sigma_p$, which is independent of the quasimomentum $\Th=(\th_1,\th_2)$, while
  the terem $p(S',\th_1,\th_2)$ determines the absolutely continuous spectrum $\sigma_{ac}$ in each case, as the quasimomentum varies in the Brillouin zone
  $[-\pi,\pi]^2$.  So we obtained the following theorem in \cite{LJL}.
 \newtheorem{th1.1}[th1.0]{Theorem}
 \begin{th1.1}
 \label{th1.1}
 Assuming all $q_i's$ are identical and even,
 \begin{enumerate}
 \item[(a)] $\displaystyle  \sigma_{ac}(H_{T})=\ \sigma_{ac}(H_{TH})=\left\{\rho^2\in \bfR:\ S'(a,\rho)\in\left[-\frac{1}{2},1\right]\right\}$.
 \item[(b)] $\displaystyle  \sigma_{ac}(H_{eT})=\left\{\rho^2\in \bfR:\ S'(a,\rho)\in\left[-\frac{3}{5},1\right]\right\}$.
 \item[(c)]   $\displaystyle   \sigma_{ac}(H_{S})=\sigma_{ac}(H_{H})=\sigma_{ac}(H_{trS})=\left\{\rho^2\in \bfR:\ S'(a,\rho)\in\left[-1,1\right]\right\}$.
 \end{enumerate}
\end{th1.1}

  In this paper, we shall study the 5 remaining Archimedean tilings.  Here the dispersion relations are even more complicated.  So we need to restrict
  $q_j$'s to be identical and even. Our method is still the characteristic function approach, putting all the information into a system of equations for the
  coefficients of the quasiperiodic solutions, and then evaluating the determinant of the resulting matrix. In all the cases, the matrices are large, even with our
  clever choice of fundamental domains.  The sizes starts with $18\times 18$, up to $36\times 36$.  So as in \cite{LJL}, we need to use the software
  Mathematica to help simplify the determinants of these big matrices.
   We first arrive at the following theorem.
 \newtheorem{th1.2}[th1.0]{Theorem}
 \begin{th1.2}
 \label{th1.2}
  Assume that all the $q_j$'s are identical (denoted as $q$), and even. We also let $\theta_1,\ \theta_2\in[-\pi,\pi]$. The dispersion relation of the
  periodic quantum graph associated with each Archimedean tiling is given by the following.
 \begin{enumerate}
 \item[(a)] For $H_{trH}$ associated with truncated hexagonal tiling ($(3,12^2)$),
  $$
  3S^3S'(3S'+2)\left\{81(S')^4-54(S')^3-45(S')^2+18S'
   -8\cos(\frac{\th_1}{2})\cos(\frac{\th_2}{2})\cos(\frac{\th_1-\th_2}{2})+8\right\}=0.
  $$
 \item[(b)] For $H_{SS}$ associated with snub square tiling ($(3,4,6,4)$), with $c=\cos \th_1,\ d=\cos\th_2$,
  \begin{eqnarray*}
  \lefteqn{S^6\,
  \left\{625 (S')^4-250 (S')^2-40 S'+1\right.}\\
 &&\left.-100(c+d)(S')^2-40(c+d+c d)S'-4(c+d+4c d-c^2-d^2)\right\}=0.
  \end{eqnarray*}
 \item[(c)] For $H_{RTH}$ associated with {rhombi}-trihexagonal tiling ($(3^2,4,3,4)$),
 \begin{eqnarray*}
  \lefteqn{S^6\, \left\{2048(S')^6-1536(S')^4-128(S')^3+192(S')^2-3-2(64(S')^3+32(S')^2-1)(\cos\th_1\right.}\\
  &&+\cos(\th_1+\th_2)+\cos\th_2)+(\cos2(\th_1+\th_2)+\cos2\th_1+\cos2\th_2)-2(\cos(2\th_1+\th_2)\\
  &&\left.+\cos(\th_1-\th_2)+\cos(\th_1+2\th_2))\right\}= 0.
   \end{eqnarray*}
 \item[(d)] For $H_{STH}$ associated with snub trihexagonal tiling ($(3^4,6)$),
 {\small
 \begin{eqnarray*}
 \lefteqn{S^9\, \left\{ 15625 (S')^6-9375 (S')^4-2000 (S')^3+675 (S')^2+120 S'-11\right.}\\
  &&+2\, (\cos2\theta_1+\cos2\theta_2+\cos2(\theta_1+\theta_2))-4\, (\cos(2\th_1+\th_2)+\cos(\th_1+2\th_2)+\cos(\th_2-\th_1))(2+5 S') \\
 &&  \left. -4\, (\cos\th_1+\cos\th_2+\cos(\th_1+\th_2))(250 (S')^3+150 (S')^2+15 S'-2)
  \right\}=0.
\end{eqnarray*}}
 \item[(e)] For $H_{trTH}$ associated with truncated trihexagonal tiling ($(4,6,12)$),
 \begin{small}
 \begin{eqnarray*}
 \lefteqn{S^6\, \left\{ 531441(S')^{12}
 -1062882(S')^{10}+728271 (S')^8-204120 (S')^6
+21627 (S')^4-918 (S')^2+15 \right.}\\
  &&+2\, (\cos2\theta_1+\cos2\theta_2+\cos2(\theta_1-\theta_2))+4\, (\cos(2\theta_1-\theta_2)+\cos(\theta_1+\theta_2)+\cos(\theta_1-2\theta_2))(1-18S'^2)
 \\
 && \left.+4\, (\cos\theta_1+\cos(\theta_1-\theta_2)+\cos\theta_2)(-2187(S')^6+1215(S')^4-135(S')^2+4)\right\}=0.
\end{eqnarray*}
\end{small}
 \end{enumerate}
 \end{th1.2}

 All the dispersion relations are of the form $S^i (3S'+2)^j\, p(S',\th_1,\th_2)=0$, where $p$ is a polynomial in $S'=S'(a,\rho)$. Thus the point spectrum of each
 periodic Schr\"odinger operator is determined by the term $S^i(3S'+2)^j$, while the absolutely continuous spectrum is determined by the part $p(S',\th_1,\th_2)$.
 As $\Th$ varies in the Brillouin zone, the range of $S'(a,\rho)$ can be derived through some tedious elementary algebra.  The following is our second main theorem.
 \newtheorem{th1.3}[th1.0]{Theorem}
 \begin{th1.3}
 \label{th1.3}
 Assuming all $q_i's$ are identical and even,
 \begin{enumerate}
 \item[(a)] $\displaystyle  \sigma_{ac}(H_{trH})=\left\{\rho^2\in \bfR:\ S'(a,\rho)\in[-\frac{2}{3},0]\bigcup\ [\frac{1}{3},1] \right\}$.
 \item[(b)] $\displaystyle  \sigma_{ac}(H_{SS})=\left\{\rho^2\in \bfR:\ S'(a,\rho)\in[-\frac{3}{5},1]\right\}$.
 \item[(c)]   $\displaystyle   \sigma_{ac}(H_{RTH})=\left\{\rho^2\in \bfR:\ S'(a,\rho)\in[-\frac{3}{4},1]\right\}$.
 \item[(d)] $\displaystyle   \sigma_{ac}(H_{STH})=\left\{\rho^2\in \bfR:\ S'(a,\rho)\in[-\frac{1+\sqrt{3}}{5},1]\right\}$.
 \item[(e)] $\displaystyle   \sigma_{ac}(H_{trTH})=\left\{\rho^2\in \bfR:\ S'(a,\rho)\in[-1,-\dfrac{1}{\sqrt{3}}]\cup [-\dfrac{1}{3},\dfrac{1}{3}]\cup[\dfrac{1}{\sqrt{3}},1]\right\}$.
 \end{enumerate}
  \end{th1.3}
There dispersion relation characterize the spectrum of each periodic Schr\"odinger operator in terms of the functions $S(a,\rho)$ and $S'(a,\rho)$,
 which are associated with the spectral problem on an interval.  In this way, the spectral problem over a quantum graph is reduced to a spectral problem
 over an interval, which we are familiar with. In particular, if $q=0$, then $\d S(a,\rho)= \frac{\sin(\rho a)}{a}$, and $S'(a,\rho)=\cos(\rho a)$.
 Thus the spectrum of each periodic operator can be computed easily (cf.\ Corollary~\ref{th7.3} below).

 We remark that these periodic quantum graphs provide a good model for the wave functions of crystal lattices. This is the so-called quantum network model (QNM)
\cite{ALM}. For example, graphene is associated with hexagonal tiling, with an identical and even potential function $\d q(x)=-0.85+\frac{d}{1.34}\, \sin^2(\frac{\pi x}{d})$,
where $d$ is the distance between neighboring atoms.  In fact under this model, several graphene-like materials, called carbon allotropes, are associated
with periodic quantum graphs \cite{K2013,EI}. Thus a study of the spectrum, which is the energy of wave functions, is important.

 In sections 2 through 6, we shall derive the dispersion relations of the periodic Schr\"odinger operators acting on the 5 above-mentioned Archimedean tilings.
 We shall invoke the vertex conditions as well as the Floquet-Bloch conditions on the boundary of the fundamental domains. With the help of Mathematica, the determinants
 the large resulting matrices are evaluated and then simplified, generating the dispersion relations. In section 7, we shall further analyze the dispersion relations to
 evaluate the point spectrum (and generate some eigenfunctions), plus the range of $S'(a,\rho)$ for each periodic Schr\"odinger operator.  In particular, we shall
 prove Theorem~\ref{th1.3}. Finally we shall have a section on concluding remarks, followed by two appendices.
 \section{Truncated hexagonal tiling}
 \vskip0.2in
  \begin{figure}[h!]
 \centering\includegraphics[width=15cm]{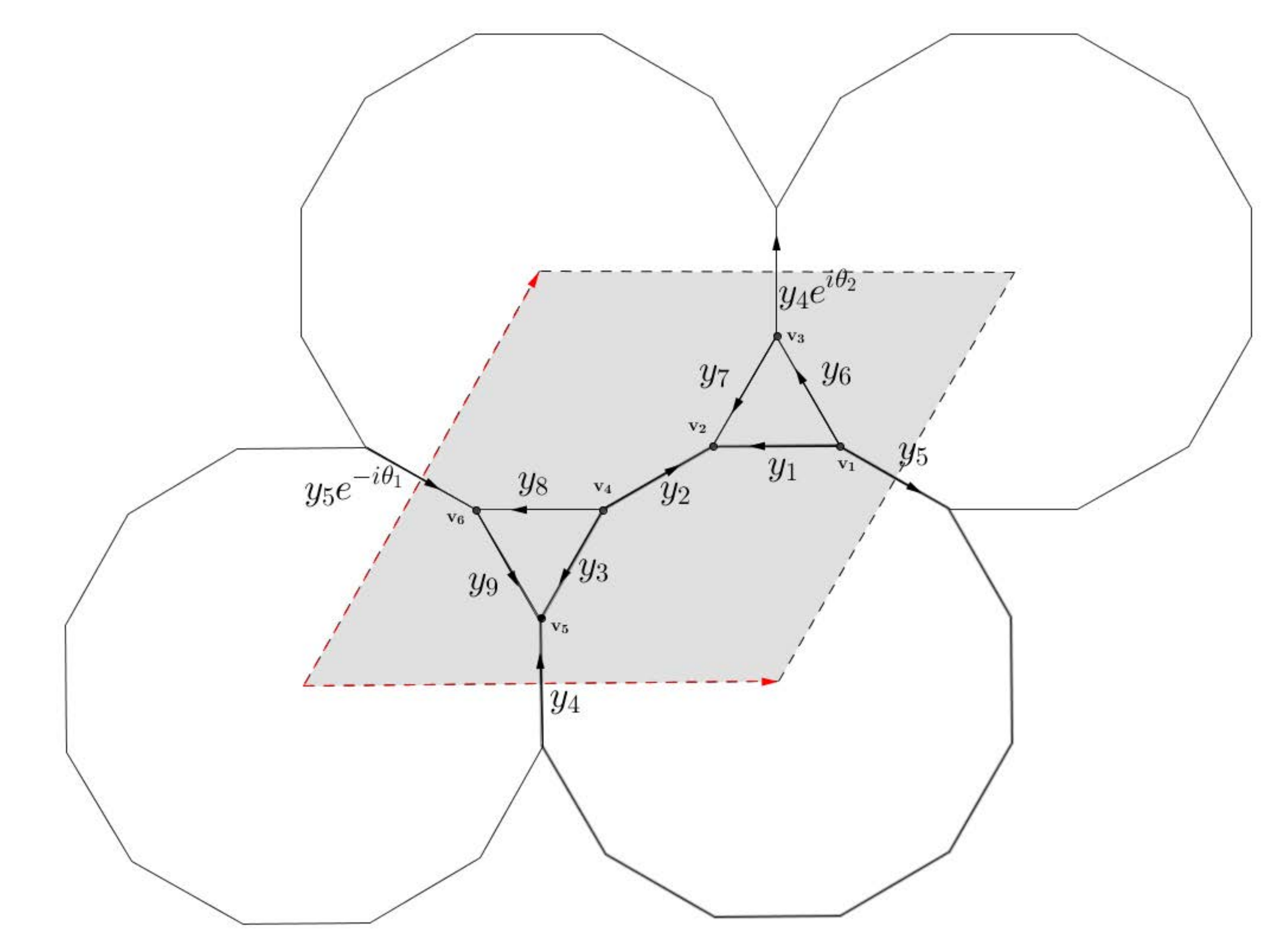}
 \caption{Fundamental domain for truncated hexagonal tiling}
 \label{fig2.1}
 \end{figure}

  The truncated hexagonal tiling $(3,12^2)$, denoted by $G$, is a periodic graph on the plane determined by one triangle and two dodecagons around a vertex.
   Each edge has length $a$.
  As shown in Fig.1, let $S$ be a parallelogram defined by two vectors $\vec{k}_1=((2+\sqrt{3})a,0)$ and $\d \vec{k}_2=((1+\frac{\sqrt{3}}{2})a,(\frac{3}{2}+\sqrt{3})a)$. Let $W=G\cap S$. Obviously the entire graph $G$ is covered by the
  translations of $W$, namely
  $$
  G=\bigcup \{ (\bp\circ W):\ \bp\in \bfZ^2\}.
  $$
 So $W$ is a fundamental domain of $G$. Applying the continuity and Kirchhoff conditions coupled with Floquet-Bloch conditions, we obtain
the following equations at:
\begin{align*}
&\textbf{v}_1:
\begin{cases}
& y_6(0)=y_1(0)=y_{5}(0); \\
& y'_6(0)+y'_1(0)+y'_{5}(0)=0.
\end{cases}\\
&\textbf{v}_2: \begin{cases}
& y_1(a)=y_{2}(a)=y_{7}(a); \\
&  y'_1(a)+y'_{2}(a)+y'_{7}(a)=0.
\end{cases}\\
&\textbf{v}_3: \begin{cases}
& y_7(0)=y_4(0)e^{i\th_2}=y_{6}(a);\\
&y'_7(0)+y'_4(0)e^{i\th_2}-y'_{6}(a)=0.
\end{cases} \\
&\textbf{v}_4: \begin{cases}
& y_2(0)=y_{3}(0)=y_{8}(0);\\
&y'_2(0)+y'_{3}(0)+y'_{8}(0)=0.
\end{cases}
\end{align*}
\newpage
 \begin{align*}
&\textbf{v}_5: \begin{cases}
& y_3(a)=y_{4}(a)=y_{9}(a);\\
& y'_3(a)+y'_{4}(a)+y'_{9}(a) =0.
\end{cases}\\
&\textbf{v}_{6}: \begin{cases}
& y_8(a)=y_5(a)e^{-i\th_1}=y_9(0); \\
& y'_8(a)+y'_5(a)e^{-i\th_1}-y'_9(0)=0.
\end{cases}
\end{align*}
Since $y_j=A_jC_j+B_jS_j$, we can rewrite the above equations as a
linear system of $(A_1,\cdots, A_{9})$ and $(B_1,\cdots, B_{9})$, thus
  \begin{equation*}
  \left\{
  \begin{array}{l}
  A_1=A_5=A_6;\\
  B_1+B_5+B_6=0;\\
  A_1C_1+B_1S_1=A_2 C_2+B_2S_2=A_7 C_7+B_7S_7;\\
  A_1 C_1'+B_1S_1'+A_2 C_2'+B_2S_2'+A_7C_7'+B_7S_7'=0;\\
  \rme^{i \th_2} A_4=A_7=A_6C_6+B_6S_6;\\
  \rme^{i \th_2} B_4+B_7-A_6C_6'-B_6S_6'=0;\\
  A_2=A_3=A_8;\\
  B_2+B_3+B_8=0;\\
  A_3C_3+B_3S_3=A_4C_4+B_4S_4=A_9C_9+B_9S_9;\\
  A_3C_3'+B_3S_3'+A_4C_4'+B_4S_4'+A_9C_9'+B_9S_9'=0;\\
  \rme^{-i\th_1}(A_5C_5+B_5S_5)=A_8C_8+B_8S_8=A_9;\\
  -B_9+\rme^{-i \th_1}(A_5 C_5'+B_5S_5')+A_8C_8'+B_8S_8'=0.
 \end{array}
 \right.
 \end{equation*}

  With $\widetilde{\al}=\rme^{-i\th_1}$ and $\be=\rme^{i\th_2}$,  the characteristic equation $\Phi(\rho)$ is given by the following determinant of a $18\times18$ matrix,
 as given in Appendix A.
 This determinant seems to be too big. So we assume all potentials $q_j$'s are identical and even. Hence, $S_j=S$, $C_j=C$ for any $j$.
 Then we use Mathematica to compute and simplify according to the coefficients $\widetilde{\al}^i\be^j$.
 It still look very wild. But after we apply the Lagrange identity $SC'=CS'-1$, things looks simpler and some pattern comes up.
 \begin{footnotesize}
 \begin{eqnarray*}
 \Phi(\rho)&=&S^3\, \left\{ -3(S'(3C+1)+C)-\be(3S'^2+S'(4+6C)+2C)+\widetilde{\al}^2\be\, (S'(6C+2)+C(3C+4))\right.\\
 &&\quad -3\widetilde{\al}^2\be^2(9S'(3C+1)+C)+\widetilde{\al} \be^2(S'(6C+2)+C(3C+4))+\widetilde{\al}(S'(6C+2)+C(3C+4))\\
 &&\quad +\widetilde{\al}\be (162 C^2 S'^4
 +81 C S'^3(5C^2-2)+S'^2(162 C^4-405 C^2-54 C+32)+S'(18+98C-54 C^2-162 C^3)\\
 &&\quad \left.+2C(16C+9)) \right\}\\
  &=& \widetilde{\al}\be S^3 \,\left\{ -6(S'(3C+1)+C)\cos(\th_2-\th_1)-2(3S'^2+(6C_4)S'+2C)\cos\th_1\right.\\
  &&\quad -2(S'(6C+2)+C(3C+4))\cos\th_2+162C^2S'^4+81 CS'^3(5C^2-2)+S'^2(162 C^4-405 C^2-54C+32)\\
  &&\left. \quad +S'(18+98C-54C^2-162C^3)+2C(16C+9)\right\}.
  \end{eqnarray*}
  \end{footnotesize}
  \noindent
  If we assume $q$ to be even so that $C=S'$, the dispersion relation for
  the periodic quantum graph associated with truncated hexagonal tiling  becomes
  \begin{footnotesize}
  \begin{eqnarray*}
  0 &=& S^3S'\{-2(9S'+6)(\cos\th_1+\cos\th_2+\cos(\th_1-\th_2))\\
  &&\quad+(9S'+6)(81(S')^4-54(S')^3-45(S')^2+18S'+6)\} \\
   &=& 3S^3S'(3S'+2)\left\{81(S')^4-54(S')^3-45(S')^2+18S'-8\cos(\dfrac{\th_1}{2})\cos(\dfrac{\th_2}{2})\cos(\dfrac{\th_1-\th_2}{2})+8\right\}.
   \end{eqnarray*}
 \end{footnotesize}

 \section{Snub square tiling}
 \vskip0.2in
   The snub square tiling $(3,4,3^2,4)$, denoted by $G$, is a periodic graph on the plane determined by one triangle and one square, then two triangles and one more square
    around a vertex. Each edge has length $a$.
  Let $S$ be an irregular hexagon as defined in Fig.2, with two translation vectors $\d \vec{k}_1=(\frac{\sqrt{3}+1}{2}a,\frac{\sqrt{3}+1}{2}a)$ and $\d \vec{k}_2=(\frac{\sqrt{3}+1}{2}a,\frac{-(\sqrt{3}+1)}{2}a)$. Let $W=G\cap S$. Obviously the entire graph $G$ is covered by the
  translations of $W$, namely
  $$
  G=\bigcup \{ (\bp\circ W):\ \bp\in \bfZ^2\}.
  $$
\begin{figure}[h!]
\centering
\includegraphics[width=10cm]{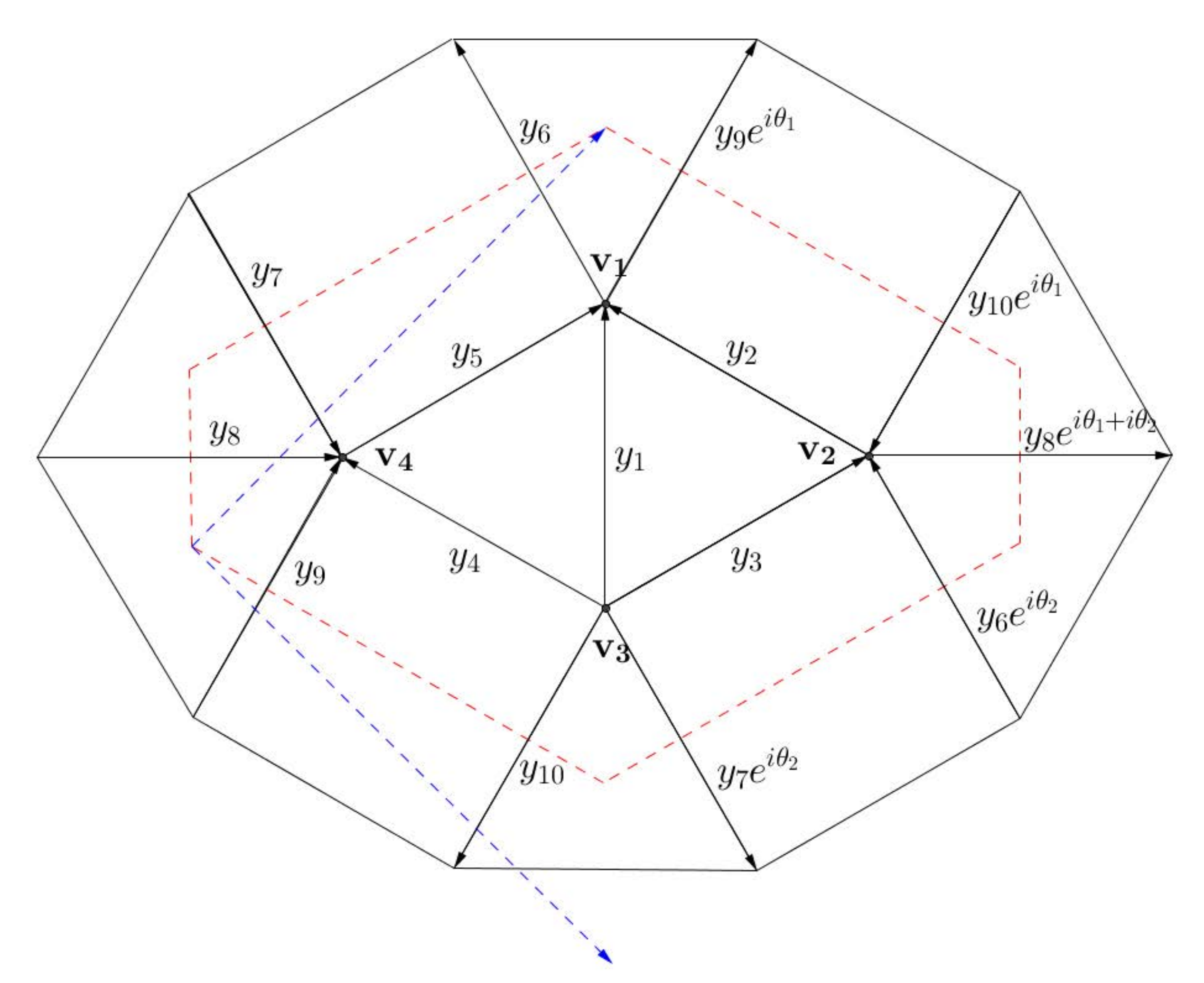}
\caption{Fundamental domain for snub square tiling}
 \label{fig3.1}
\end{figure}
  So $W$ is a fundamental domain of $G$. The continuity and Kirchhoff conditions together with the Floquet-Bloch conditions, yield the following equations at:
  \begin{align*}
  &\textbf{v}_1:
  \begin{cases}
  & y_1(a)=y_2(a)=y_5(a)=y_6(0)=y_9(0)e^{i\th_1}; \\
  & y'_1(a)+y'_2(a)+y'_5(a)-y'_6(0)-y'_9(0)e^{i\th_1}=0.
  \end{cases}\\
  &\textbf{v}_2: \begin{cases}
  & y_2(0)=y_3(a)=y_6(a)e^{i\th_2}=y_8(0)e^{i\th_1+i\th_2}=y_{10}(a)e^{i\th_1}; \\
  &  -y'_2(0)+y'_3(a)+y'_6(a)e^{i\th_2}-y'_8(0)e^{i\th_1+i\th_2}+y'_{10}(a)e^{i\th_1}=0.
  \end{cases}\\
  &\textbf{v}_3: \begin{cases}
  & y_1(0)=y_3(0)=y_4(0)=y_7(0)e^{i\th_2}=y_{10}(0);\\
  &y'_1(0)+y'_3(0)+y'_4(0)+y'_7(0)e^{i\th_2}+y'_{10}(0)=0.
  \end{cases} \\
  &\textbf{v}_4: \begin{cases}
  & y_4(a)=y_5(0)=y_7(a)=y_8(a)=y_9(a);\\
  &y'_4(a)-y'_5(0)+y'_7(a)+y'_8(a)+y'_9(a)=0.
  \end{cases}
  \end{align*}
  Since $y_j=A_jC_j+B_jS_j$, we can rewrite the above equations as a
  	linear system of $(A_1,\cdots, A_{10})$ and $(B_1,\cdots, B_{10})$, thus
\[\begin{cases}
&A_1C_1+B_1S_1=A_2C_2+B_2S_2=A_5C_5+B_5S_5=A_6=A_9e^{i\th_1};\\
&A_1C_1+B_1S_1+A_2C'_2+B_2S'_2+A_5C'_5+B_5S'_5-B_6-B_9e^{i\th_1}=0;\\
&A_3C_3+B_3S_3=A_2=(A_{10}C_{10}+B_{10}S_{10})e^{i\th_1}=A_8e^{i(\th_1+\th_2)}=(A_6C_6+B_6S_6)e^{i\th_2};\\
&A_3C'_3+B_3S'_3-B_2+(A_{10}C'_{10}+B_{10}S'_{10})e^{i\th_1}-B_8e^{i(\th_1+\th_2)}+(A_6C'_6+B_6S'_6)e^{i\th_2}=0;\\

&A_1=A_3=A_4=A_{10}=A_7e^{i\th_2};\\
&B_1+B_3+B_4+B_{10}+B_7e^{i\th_2}=0\\
&A_5=A_7C_7+B_7S_7=A_8C_8+B_8S_8=A_9C_9+B_9S_9=A_4C_4+B_4S_4;\\
&-B_5+A_7C_7+B_7S_7+A_8C_8+B_8S_8+A_9C_9+B_9S_9+A_4C_4+B_4S_4=0.
\end{cases}\]
Let $\al=e^{i\th_1}$ and $\be=e^{i\th_2}$. Then the characteristic function $\Phi(\rho)$ is given by the determinant of a $20\times 20$ matrix.

This determinant is quite big to handle. We simplify it by making the $S_j's$ and $C_j's$ to be the same as we did in the previous section. Then we use Mathematica to compute and simplify according to the coefficients $\al^i\be^i$. After the computation, it still looks complicated, involving around 300 terms. However, by the help of the Lagrange identity, $SC'=CS'-1$, things looks  a lot simpler and some pattern comes out. The characteristic function will be
\begin{footnotesize}
\begin{align*}
\Phi(\rho)=&(\al^5 \be^3+\al \be^3) S^6+(\al^3 \be^5+\al^3 \be) S^6-(\al^4 \be^3+\al^2 \be^3) S^6 \left(12 \left(S'\right)^2+2 (13 C+5) S' +2 \left(6 C^2+5 C+1\right)\right)\\
&-(\al^4 \be^4+\al^2 \be^2) S^6 \left(3 S'+7 C+4\right)-(\al^4 \be^2+\al^2 \be^4) S^6 \left(7 S'+3 C+4\right)\\
&-(\al^3 \be^4+\al^3 \be^2) S^6 \left(9 \left(S'\right)^2+(32 C+10) S'+9 C^2+10 C+2\right)\\
&+\al^3 \be^3 S^6 \left(180 C \left(S'\right)^3-3 \left(18-95 C^2\right) \left(S'\right)^2+2 \left(70 C^3-71 C-10\right) S'+20 C^4-54 C^2-20 C+5\right)\\
=&\al^3\be^3S^6\left(2\cos2\th_1+2\cos2\th_2-2\cos\th_1\left(12 \left(S'\right)^2+2 (13 C+5) S' +2 \left(6 C^2+5 C+1\right)\right)\right.\\
&\left.-2\cos(\th_1+\th_2)\left(3 S'+7 C+4\right)-2\cos(\th_1-\th_2)\left(7 S'+3 C+4\right)\right.\\
&\left.-2\cos\th_2\left(9 \left(S'\right)^2+(32 C+10) S'+9 C^2+10 C+2\right)\right.\\
&\left.+\left(180 C \left(S'\right)^3-3 \left(18-95 C^2\right) \left(S'\right)^2+2 \left(70 C^3-71 C-10\right) S'+20 C^4-54 C^2-20 C+5\right)\right).
\end{align*}
\end{footnotesize}
Now, if we assume $q$ to be even so that $C=S'$, the dispersion relation for the periodic quantum graph associated with snub square tiling will be
\begin{align*}
0=&S^6\left(2\cos2\th_1+2\cos2\th_2-2(\cos\th_1+\cos\th_2)(50(S')^2+20S'+2)\right.\\
&\left.-2(\cos(\th_1+\th_2)+\cos(\th_1-\th_2))(10S'+4)+625 (S')^4-250 (S')^2-40 S'+5\right).
\end{align*}

 Let $c=\cos\th_1$ and $d=\cos\th_2$. This yields the dispersion relation in Theorem~\ref{th1.2}(b).

 \section{Rhombi-trihexagonal tiling}
\vskip0.2in The rhombi-trihexagonal tiling $(3,4,6,4)$, denoted by $G$, is a periodic graph on the plane determined by one triangle, one square, one hexagon and
  another square around a vertex.
   Each edge has length $a$.
  As shown in fig.3, let $S$ be a parallelogram defined by two vectors $\vec{k}_1=\left(\dfrac{\sqrt{3}+1}{2}a,\dfrac{\sqrt{3}+3}{2}a\right)$ and $\vec{k}_2=\left(\dfrac{\sqrt{3}+1}{2}a,-\dfrac{\sqrt{3}+3}{2}a\right)$. Let $W=G\cap S$. Obviously the entire graph $G$ is covered by the
  translations of $W$, namely
  $$
  G=\bigcup \{ (\bp\circ W):\ \bp\in \bfZ^2\}.
  $$
\begin{figure}[h!]
\centering
\includegraphics[width=16cm]{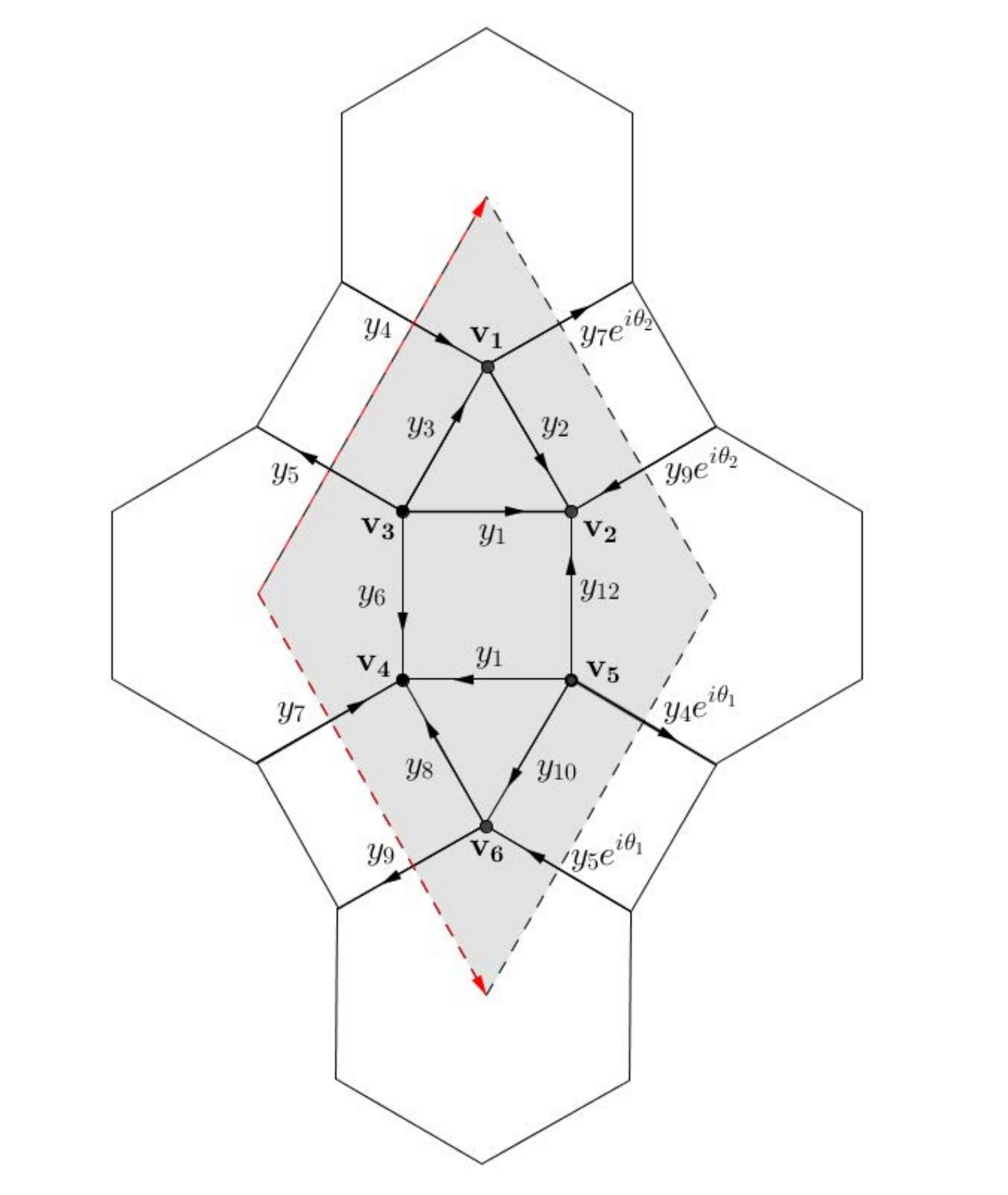}
\caption{Fundamental domain of rhombi-trihexagonal tiling}
\label{fig4.1}
\end{figure}
 This $W$ is a fundamental domain of $G$. The continuity and Kirchhoff conditions coupled with Floquet-Bloch conditions yield:
\begin{align*}
&\textbf{v}_1:
\begin{cases}
& y_2(0)=y_3(a)=y_{4}(a)=y_{7}(0)e^{i\th_2}; \\
& -y'_2(0)+y'_3(a)+y'_{4}(a)-y'_{7}(0)e^{i\th_2}=0.
\end{cases}\\
&\textbf{v}_2: \begin{cases}
& y_1(a)=y_{2}(a)=y_{12}(a)=y_9(a)e^{i\th_2}; \\
& y'_1(a)+y'_{2}(a)+y'_{12}(a)+y'_9(a)e^{i\th_2}=0.
\end{cases}\\
&\textbf{v}_3: \begin{cases}
& y_1(0)=y_3(0)=y_{5}(0)=y_{6}(0);\\
&y'_1(0)+y'_3(0)+y'_{5}(0)+y'_{6}(0)=0.
\end{cases} \\
&\textbf{v}_4: \begin{cases}
& y_6(a)=y_{7}(a)=y_{8}(a)=y_{11}(a); \\
&y'_6(a)+y'_{7}(a)+y'_{8}(a)+y'_{11}(a)=0.
\end{cases} \\
&\textbf{v}_5: \begin{cases}
& y_4(0)e^{i\th_1}=y_{10}(0)=y_{11}(0)=y_{12}(0);\\
& y'_4(0)e^{i\th_1}+y'_{10}(0)+y'_{11}(0)+y'_{12}(0) =0.
\end{cases} \\
&\textbf{v}_{6}: \begin{cases}
& y_8(0)=y_9(0)=y_{10}(a)=y_5(a)e^{i\th_1}; \\
& -y'_8(0)-y'_9(0)+y'_{10}(a)+y'_5(a)e^{i\th_1}=0.
\end{cases}&
\end{align*}
Since $y_j=A_jC_j+B_jS_j$, we can rewrite the above equations as a
linear system of $(A_1,\cdots, A_{18})$ and $(B_1,\cdots, B_{18})$ with $%
\alpha=e^{i\theta_1}$ and $\beta=e^{i\theta_2}$. Thus,
\begin{equation}
\begin{cases}
&A_2=A_3C_3+B_3S_3=A_4C_4+B_4S_4=A_7\beta;\\
&-B_2+A_3C'_3+B_3S'_3+A_4C'_4+B_4S'_4-B_7\beta=0;\\
&A_1C_1+B_1S_1=A_2C_2+B_2S_2=A_{12}C_{12}+B_{12}S_{12}=(A_9C_9+B_9S_9)\beta;\\
&A_1C'_1+B_1S'_1+A_2C'_2+B_2S'_2+A_{12}C'_{12}+B_{12}S'_{12}+(A_9C'_9+B_9S'_9)\beta=0;\\
&A_1=A_3=A_5=A_6;\\
&B_1+B_3+B_5+B_6=0;\\
&A_6C_6+B_6S_6=A_7C_7+B_7S_7=A_8C_8+B_8S_8=A_{11}C_{11}+B_{11}S_{11};\\
&A_6C'_6+B_6S'_6+A_7C'_7+B_7S'_7+A_8C'_8+B_8S'_8+A_{11}C'_{11}+B_{11}S'_{11}=0;\\
&A_4\alpha=A_{10}=A_{11}=A_{12};\\
&B_4\alpha+B_{10}+B_{11}+B_{12}=0;\\
&A_8=A_9=A_{10}C_{10}+B_{10}S_{10}=(A_5C_5+B_5S_5)\alpha;\\
&-B_8-B_9+A_{10}C'_{10}+B_{10}S'_{10}+(A_5C'_5+B_5S'_5)\alpha=0.
\end{cases}
\end{equation}
Assume that the potential $q_i's$ are identical, that is, $S:=S_i$ and $C:=C_i$.
So the characteristic equation $\Phi(\rho)$ is a determinant of $24\times 24$
matrix.
With the help of Mathematica and  Lagrange identity, $SC'=CS'-1$, we evaluate the characteristics function as
{\footnotesize
  \setlength{\abovedisplayskip}{6pt}
  \setlength{\belowdisplayskip}{\abovedisplayskip}
  \setlength{\abovedisplayshortskip}{0pt}
  \setlength{\belowdisplayshortskip}{3pt}
\begin{align*}
\Phi(\rho)=&S^6\{-(1+\al^4\be^4)-(\be^2+\al^4\be^2)-(\al^2+\al^2\be^4)+2(\be+\al^4\be^3)+2(\al+\al^3\be^4)+2(\al\be^3+\al^3\be)\\
&+2(\al^3\be^3+\al\be)(32C(S')^2+32C(1+C)S'-1)\\
&+2(\al\be^2+\al^3\be^2)(16(1+2C)(S')^2+32C^2S'+16CS'-1)\\
&+2(\al^2\be+\al^2\be^3)(21C(S')^2+16C(1+2C)S'+16CS'-1)\\
&+\al^2\be^2(-1024C^2(S')^4-256C(-3+8C^2)(S')^3+32(-24C+48C^2-32C^4)(S')^2\\
&+128C(-2+C+6C^2)S'-64C^2+6)\}\\
=&\al^2\be^2S^6\{-(\al^{-2}\be^{-2}+\al^2\be^2)-(\be^{-2}+\be^2)-(\al^{-2}+\al^2)+2((\al^{-1}\be^{-2}+\al\be^2)+(\al^{-1}\be+\al\be^{-1})\\
&+2(\al^{-2}\be^{-1}+\al^2\be))+2(\al+\al^{-1})(32C(S')^2+32C(1+C)S'-1)\\
&+2(\al^{-1}\be^{-1}+\al\be)(16(1+2C)(S')^2+32C^2S'+16CS'-1)\\
&+2(\be+\be^{-1})(21C(S')^2+16C(1+2C)S'+16CS'-1)\\
&-1024C^2(S')^4-256C(-3+8C^2)(S')^3+32(-24C+48C^2-32C^4)(S')^2\\
&+128C(-2+C+6C^2)S'-64C^2+6\}\\
=&\al^2\be^2S^6\{-2(\cos2(\th_1+\th_2)+\cos2\th_1+\cos2\th_2)+4(\cos(2\th_1+\th_2)+\cos(\th_1-\th_2)+\cos(\th_1+2\th_2))\\
&+4\cos\th_1(32C(S')^2+32C(1+C)S'-1)+4\cos(\th_1+\th_1)(16(1+2C)(S')^2+32C^2S'+16CS'-1)\\
&+4\cos\th_2(21C(S')^2+16C(1+2C)S'+16CS'-1)\\
&-1024C^2(S')^4-256C(-3+8C^2)(S')^3+32(-24C+48C^2-32C^4)(S')^2\\
&+128C(-2+C+6C^2)S'-64C^2+6\}.
\end{align*}
}
\indent
If we further assume that the the potential $q$  is even then $C=S'$, and it yields the dispersion relation in
  Theorem \ref{th1.2}(c).

 \section{Snub trihexagonal tiling}
 \hskip0.25in
 The snub trihexagonal tiling $(3^3,6)$, denoted by $G$, is a periodic graph on the plane determined by three triangles and one hexagon
    around a vertex. Each edge has length $a$.
  Let $S$ be a regular hexagon as defined in Fig.4, with two translation vectors $\vec{k}_1=(2a,\sqrt{3}a)$ and $\vec{k}_2=(\frac{5a}{2},-\frac{\sqrt{3}a}{2})$. Let $W=G\cap S$ is a fundamental domain of $G$.
\begin{figure}[h!]
\centering
\includegraphics[width=16cm]{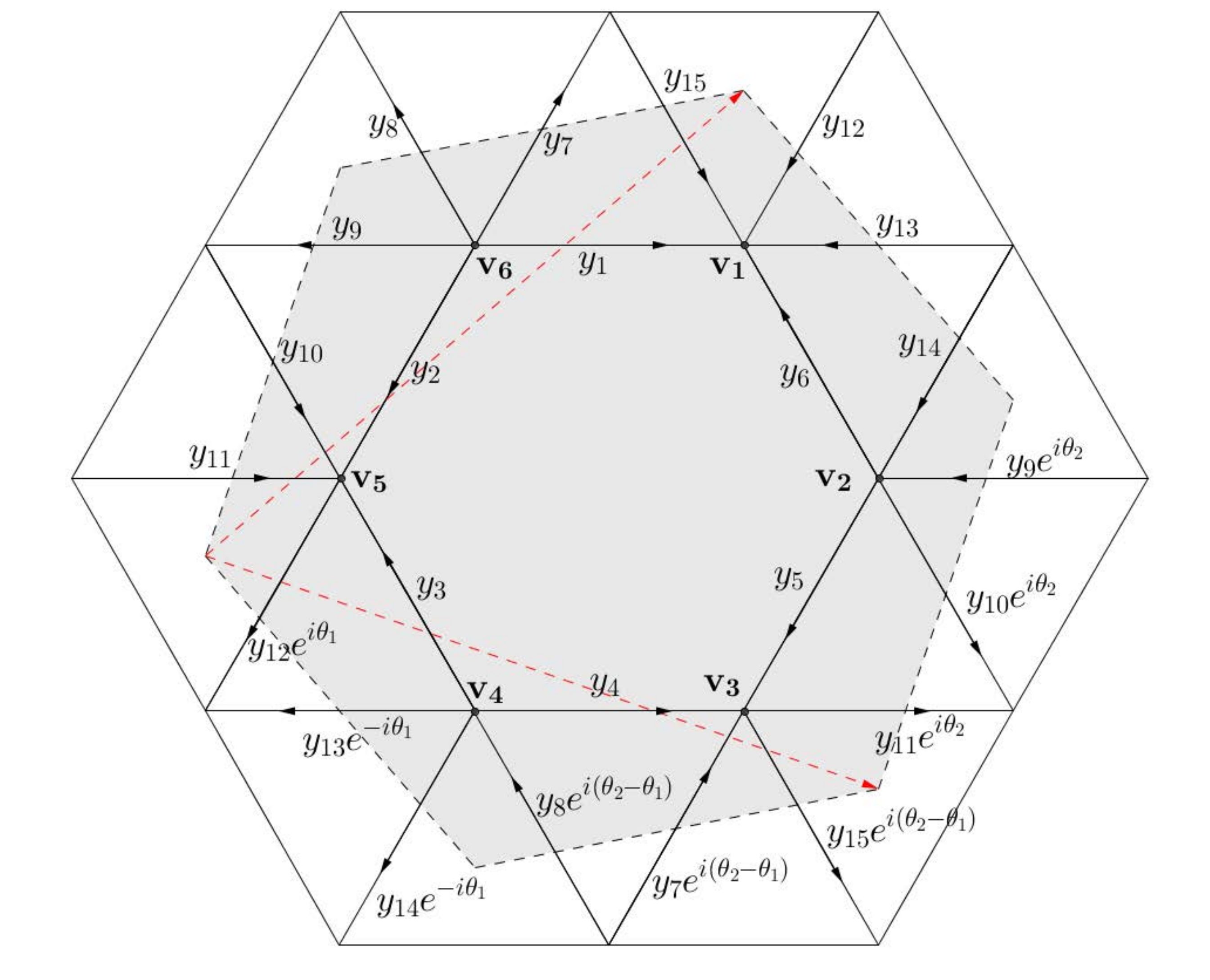}
\caption{Fundamental domain of snub trihexagonal tiling}
\label{fig5.1}
\end{figure}
 We apply the continuity and Kirchhoff conditions coupled with Floquet-Bloch conditions to obtain:
\begin{align*}
&\textbf{v}_1:
\begin{cases}
& y_1(a)=y_6(a)=y_{15}(a)=y_{12}(a)=y_{13}(a); \\
& y'_1(a)+y'_6(a)+y'_{15}(a)+y'_{12}(a)+y'_{13}(a)=0.
\end{cases}\\
&\textbf{v}_2: \begin{cases}
& y_6(0)=y_5(0)=y_{10}(0)e^{i\th_2}=y_{14}(a)=y_9(a)e^{i\th_2}; \\
& y'_6(0)+y'_5(0)+y'_{10}(0)e^{i\th_2}-y'_{14}(a)-y_9(a)e^{i\th_2} =0.
\end{cases}\\
&\textbf{v}_3: \begin{cases}
& y_4(a)=y_5(a)=y_7(a)e^{i(\th_2-\th_1)}=y_{15}(0)e^{i(\th_2-\th_1)}=y_{11}(0)e^{i\th_2};\\
& y'_4(a)+y'_5(a)+y'_7(a)e^{i(\th_2-\th_1)}-y'_{15}(0)e^{i(\th_2-\th_1)}-y'_{11}(0)e^{i\th_2} =0.
\end{cases}
 \end{align*}
 \newpage
 \begin{align*}
 &\textbf{v}_4: \begin{cases}
 & y_3(0)=y_4(0)=y_{13}(0)e^{-i\th_1}=y_{14}(0)e^{-i\th_1}=y_8(a)e^{i(\th_2-\th_1)}; \\
 &y'_3(0)+y'_4(0)+y'_{13}e^{-i\th_1}+y'_{14}(0)e^{-i\th_1}-y'_8(a)e^{i(\th_2-\th_1)} =0.
 \end{cases} \\
 &\textbf{v}_5: \begin{cases}
 & y_2(a)=y_3(a)=y_{10}(a)=y_{11}(a)=y_{12}(0)e^{-i\th_1};\\
 &  y'_2(a)+y'_3(a)+y'_{10}(a)+y'_{11}(a)-y'_{12}(0)e^{-i\th_1} =0.
 \end{cases}\\
&\textbf{v}_{6}: \begin{cases}
& y_1(0)=y_2(0)=y_7(0)=y_8(0)=y_9(0); \\
& y'_1(0)+y'_2(0)+y'_7(0)+y'_8(0)+y'_9(0)=0.
\end{cases}&
\end{align*}
Since $y_j=A_jC_j+B_jS_j$, we can rewrite the above equations as a
linear system of $(A_1,\cdots, A_{18})$ and $(B_1,\cdots, B_{18})$ with $%
\widetilde{\al}=e^{-i\theta_1}$ and $\beta=e^{i\theta_2}$. Thus,
\begin{equation}
\begin{cases}
& A_1C_1+B_1S_1=A_6C_6+B_6S_6=A_{15}C_{15}+B_{15}S_{15}=A_{12}C_{12}+B_{12}S_{12}=A_{13}C_{13}+B_{13}S_{13};\\
&A_1C'_1+B_1S'_1+A_6C'_6+B_6S'_6+A_{15}C'_{15}+B_{15}S'_{15}+A_{12}C'_{12}+B_{12}S'_{12}+A_{13}C'_{13}+B_{13}S'_{13}=0;\\
& A_6=A_5=\beta A_{10}=A_{14}C_{14}+B_{14}S_{14}=(A_{9}C_{9}+\beta B_{9}S_{9});\\
&B_6+B_5+\beta B_{10}-(A_{14}C'_{14}+B_{14}S'_{14})-\beta (A_{9}C'_{9}+B_{9}S'_{9})=0; \\
&A_4C_4+B_4S_4=A_5C_5+B_5S_5=\widetilde{\al}\beta(A_7C_7+B_7S_7)=\widetilde{\al}\beta A_{15}=\beta A_{11};\\
&A_4C'_4+B_4S'_4+A_5C'_5+B_5S'_5+\widetilde{\al}\beta(A_7C'_7+B_7S'_7)-\widetilde{\al}\beta B_{15}-\beta B_{11}=0;\\
&A_3=A_4=\widetilde{\al} A_{13}=\widetilde{\al}A_{14}=\widetilde{\al}\beta(A_8C_8+B_8S_8);\\
&B_3+B_4+\widetilde{\al} B_{13}+\widetilde{\al}B_{14}-\widetilde{\al}\beta(A_8C'_8+B_8S'_8)=0;\\
&A_2C_2+B_2S_2=A_3C_3+B_3S_3=A_{10}C_{10}+B_{10}S_{10}=A_{11}C_{11}+B_{11}S_{11}=\widetilde{\al}A_{12};\\
&A_2C'_2+B_2S'_2+A_3C'_3+B_3S'_3+A_{10}C'_{10}+B_{10}S'_{10}+A_{11}C'_{11}+B_{11}S'_{11}-\widetilde{\al}B_{12}=0;\\
&A_1=A_2=A_7=A_8=A_9;\\
&B_1+B_2+B_7+B_8+B_9=0.
\end{cases}
\end{equation}
Next we assume that the potential $q_i's$ are identical so that $S_i=S$ and $C_i=C$. The characteristic equation $\Phi(\rho)$ is thus a determinant of $30\times 30$
matrix.
With the help of Mathematica and Lagrange identity, we evaluate the characteristic function as
{\footnotesize
  \setlength{\abovedisplayskip}{6pt}
  \setlength{\belowdisplayskip}{\abovedisplayskip}
  \setlength{\abovedisplayshortskip}{0pt}
  \setlength{\belowdisplayshortskip}{3pt}
\begin{align*}
\Phi(\rho)=&S^9\{((\widetilde{\al}^4 \be^4+\widetilde{\al}^8\be^8)+(\widetilde{\al}^4 \be^6+\widetilde{\al}^8\be^6)+(\widetilde{\al}^6 \be^4+\widetilde{\al}^6 \be^8))+((\widetilde{\al}^4 \be^5+\widetilde{\al}^8 \be^7)\\
&+(\widetilde{\al}^7 \be^5+\widetilde{\al}^5 \be^7))(-4 (C+1)-6 S')+(\widetilde{\al}^7 \be^8+\widetilde{\al}^5 \be^4) (-4-7 C-3 S')+(\widetilde{\al}^7 \be^6+\widetilde{\al}^5 \be^6) \\
& (4-18 C-47 C^2-24 C^3+(-163 C^2-176 C-12) S'-(77+277 C)(S')^2-36 (S')^3)\\
&+(\widetilde{\al}^7 \be^7+\widetilde{\al}^5 \be^5)  (4-15 C-83 C^2-60 C^3+(-235 C^2-164 C-15) S'\\
&-(175 C+53)(S')^2-30 (S')^3)+(\widetilde{\al}^6 \be^5+\widetilde{\al}^6 \be^7)  (4-12 C-63 C^2-29 C^3\\
&-6 (38 C^2+29 C+3) S'-9 (23 C+7) (S')^2-36 (S')^3)+\widetilde{\al}^6 \be^6 (-11+60 C+131 C^2\\
&-224 C^3-274 C^4+(600 C^5-2279 C^3-798 C^2+413 C+60) S'+(3850 C^4-4269 C^2\\
&-732 C+131) (S')^2+(6725 C^3-2279 C-246) (S')^3+2 (1925 C^2-137) (S')^4+600 C (S')^5)\}\\
=&\widetilde{\al}^6\be^6S^9\{((\widetilde{\al}^{-2} \be^{-2}+\widetilde{\al}^{2} \be^{2})+(\widetilde{\al}^{-2}+\widetilde{\al}^2)+(\be^{-2}+\be^2))+((\widetilde{\al}^{-2}\be^{-1}+\widetilde{\al}^2\be)\\
&+(\widetilde{\al}\be^{-1}+\widetilde{\al}^{-1}\be))(-4 (C+1)-6 S')+(\widetilde{\al}\be^2+\widetilde{\al}^{-1}\be^{-2}) (-4-7 C-3 S')+(\widetilde{\al}+\widetilde{\al}^{-1})\\
& (4-18 C-47 C^2-24 C^3+(-163 C^2-176 C-12) S'-(77+277 C)(S')^2-36 (S')^3)\\
&+(\widetilde{\al}\be+(\widetilde{\al}\be)^{-1}) (4-15 C-83 C^2-60 C^3+(-235 C^2-164 C-15) S'+(-175 C-53)(S')^2\\
&-30 (S')^3)+(\be^{-1}+\be) (4-12 C-63 C^2-29 C^3-6 (38 C^2+29 C+3) S'-9 (23 C+7) (S')^2\\
&-36 (S')^3)-11+60 C+131 C^2-224 C^3-274 C^4+(600 C^5-2279 C^3-798 C^2+413 C+60) S'\\
&+(3850 C^4-4269 C^2-732 C+131) (S')^2+(6725 C^3-2279 C-246) (S')^3\\
&+2 (1925 C^2-137) (S')^4+600 C (S')^5\}\\
=&\widetilde{\al}^6\be^6S^9\{(2\cos2(\theta_1+\theta_2)+2\cos2\theta_1+(2\cos2\theta_2))+((2\cos(2\th_1+\th_2))\\
&+2\cos(\th_2-\th_1))(-4 (C+1)-6 S')+2\cos(\th_1+2\th_2)(-4-7 C-3 S')\\
&+2\cos\th_1(4-18 C-47 C^2-24 C^3+(-163 C^2-176 C-12) S'-(77+277 C)(S')^2-36 (S')^3)\\
&+2\cos(\th_1+\th_2)(4-15 C-83 C^2-60 C^3+(-235 C^2-164 C-15) S'+(-175 C-53)(S')^2-30 (S')^3)\\
&+2\cos\th_2(4-12 C-63 C^2-29 C^3-6 (38 C^2+29 C+3) S'-9 (23 C+7) (S')^2-36 (S')^3)\\
&-11+60 C+131 C^2-224 C^3-274 C^4+(600 C^5-2279 C^3-798 C^2+413 C+60) S'\\
&+(3850 C^4-4269 C^2-732 C+131) (S')^2+(6725 C^3-2279 C-246) (S')^3\\
&+2 (1925 C^2-137) (S')^4+600 C (S')^5\}\ .
\end{align*}
}
\indent
If the potential $q$  is even then $C=S'$. Therefore
 Theorem \ref{th1.2}(d) is valid.

\section{Truncated trihexagonal tiling}
 \vskip0.2in
 The truncated trihexagonal tiling $(4,6,12)$, denoted by $G$, is a periodic graph on the plane determined by one square, one hexagon and one dodecagon
    around a vertex. Each edge has length $a$.
  As shown in fig.4, let $S$ be a parallelogram defined by two vectors $\vec{k}_1=\left(\dfrac{3+3\sqrt{3}}{2}a,\dfrac{3+\sqrt{3}}{2}a\right)$ and $\vec{k}_2=\left(\dfrac{3+3\sqrt{3}}{2}a,-\dfrac{3+\sqrt{3}}{2}a\right)$. Let $W=G\cap S$ is a fundamental domain of $G$.
\begin{figure}[h!]
\centering
\includegraphics[width=18cm]{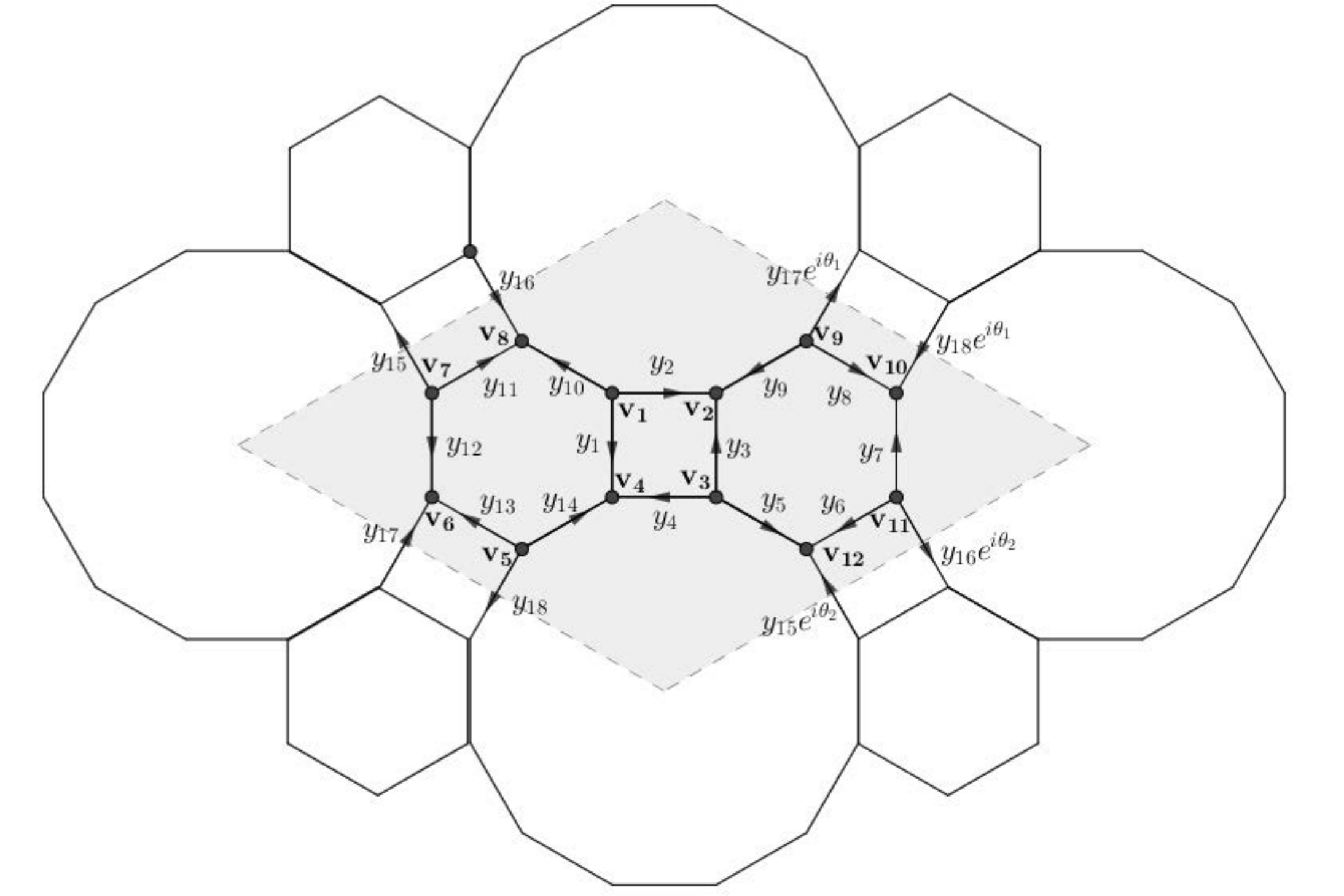}
\caption{Fundamental domain of truncated trihexagonal tiling}
\label{fig6.1}
\end{figure}
Applying the continuity and Kirchhoff conditions coupled with Floquet-Bloch conditions, we obtain:
\begin{align*}
&\textbf{v}_1:
\begin{cases}
& y_1(0)=y_2(0)=y_{10}(0); \\
& y'_1(0)+y'_2(0)+y'_{10}(0) =0.
\end{cases}& \textbf{v}_2:
\begin{cases}
& y_{2}(a)=y_{3}(a)=y_{9}(a); \\
& y'_{2}(a)+y'_{3}(a)+y'_{9}(a) =0.
\end{cases}\\
&\textbf{v}_3: \begin{cases}
& y_3(0)=y_4(0)=y_{5}(0); \\
& y'_3(0)+y'_4(0)+y'_{5}(0) =0.
\end{cases}&\textbf{v}_4:
\begin{cases}
& y_{1}(a)=y_{4}(a)=y_{14}(a) \\
& y'_{1}(a)+y'_{4}(a)+y'_{14}(a) =0
\end{cases}\\
&\textbf{v}_5: \begin{cases}
& y_{18}(0)=y_{13}(0)=y_{14}(0); \\
& y'_{18}(0)+y'_{13}(0)+y'_{14}(0) =0.
\end{cases} &\textbf{v}_6:
\begin{cases}
& y_{12}(a)=y_{13}(a)=y_{17}(a); \\
& y'_{12}(a)+y'_{13}(a)+y'_{17}(a) =0.
\end{cases}\\
&\textbf{v}_7: \begin{cases}
& y_{15}(0)=y_{12}(0)=y_{11}(0); \\
& y'_{15}(0)+y'_{12}(0)+y'_{11}(0) =0.
\end{cases} &\textbf{v}_8:
\begin{cases}
& y_{16}(a)=y_{11}(a)=y_{10}(a); \\
& y'_{16}(a)+y'_{11}(a)+y'_{10}(a) =0.
\end{cases}\\
&\textbf{v}_9: \begin{cases}
& y_{9}(0)=y_{17}(0)e^{i\theta_1}=y_{8}(0); \\
& y'_{9}(0)+y'_{17}(0)e^{i\theta_1}+y'_{8}(0) =0.
\end{cases} &\textbf{v}_{10}:
\begin{cases}
& y_{8}(a)=y_{7}(a)=y_{18}(a)e^{i\theta_1}; \\
& y'_{8}(a)+y'_{7}(a)+y'_{18}(a)e^{i\theta_1} =0.
\end{cases} \\
 &\textbf{v}_{11}: \begin{cases}
& y_{7}(0)=y_{16}(0)e^{i\theta_2}=y_{6}(0); \\
& y'_{7}(0)+y'_{16}(0)e^{i\theta_2}+y'_{6}(0) =0.
\end{cases}& \textbf{v}_{12}:
\begin{cases}
& y_{5}(a)=y_{6}(a)=y_{15}(a)e^{i\theta_2}; \\
& y'_{5}(a)+y'_{6}(a)+y'_{15}(a) e^{i\theta_2}=0.
\end{cases}
\end{align*}
 \vskip0.3in
Since $y_j=A_jC_j+B_jS_j$, we rewrite the above equations as a
linear system of $(A_1,\cdots, A_{18})$ and $(B_1,\cdots, B_{18})$ with $%
\al=e^{i\theta_1}$ and $\beta=e^{i\theta_2}$. Thus,
\begin{equation}
\begin{cases}
& A_{1}=A_{2}=A_{10}\text{; }A_{3}=A_{4}=A_{5}\text{; }A_{18}=A_{13}=A_{14};\\
&A_{15}=A_{12}=A_{11}\text{; }A_{9}=\alpha A_{17}=A_{8}\text{; }A_{7}=A_{6}=\beta A_{16}; \\
& B_{1}+B_{2}+B_{10}=0\text{; }B_{3}+B_{4}+B_{5}=0\text{; }B_{18}+B_{13}+B_{14}=0;\\
&B_{15}+B_{12}+B_{11}=0\text{; }B_{9}+\alpha B_{17}+B_{8}=0\text{; }B_{7}+B_{6}+\beta B_{16}=0; \\
&A_2C_2+B_2S_2=A_3C_3+B_3S_3=A_9C_9+B_9S_9;\\
&A_{1}C_{1}+B_{1}S_{1}=A_{4}C_{4}+B_{4}S_{4}=A_{14}C_{14}+B_{14}S_{14};\\
&A_{12}C_{12}+B_{12}S_{12}=A_{13}C_{13}+B_{13}S_{13}=A_{17}C_{17}+B_{17}S_{17};\\
&A_{16}C_{16}+B_{16}S_{16}=A_{11}C_{11}+B_{11}S_{11}=A_{10}C_{10}+B_{10}S_{10};\\
&A_{8}C_{8}+B_{8}S_{8}=A_{7}C_{7}+B_{7}S_{7}=\al(A_{18}C_{18}+B_{18}S_{18});\\
&A_{5}C_{5}+B_{5}S_{5}=A_{6}C_{6}+B_{6}S_{6}=\beta(A_{15}C_{15}+B_{15}S_{15});\\
&A_2C'_2+B_2S'_2+A_3C'_3+B_3S'_3+A_9C'_9+B_9S'_9=0;\\
&A_{1}C'_{1}+B_{1}S'_{1}+A_{4}C'_{4}+B_{4}S'_{4}+A_{14}C'_{14}+B_{14}S'_{14}=0;\\
&A_{12}C'_{12}+B_{12}S'_{12}+A_{13}C'_{13}+B_{13}S'_{13}+A_{17}C'_{17}+B_{17}S'_{17}=0;\\
&A_{16}C'_{16}+B_{16}S'_{16}+A_{11}C'_{11}+B_{11}S'_{11}+A_{10}C'_{10}+B_{10}S'_{10}=0;\\
&A_{8}C'_{8}+B_{8}S'_{8}+A_{7}C'_{7}+B_{7}S'_{7}+(A_{18}C'_{18}+\al B_{18}S'_{18})=0;\\
&A_{5}C'_{5}+B_{5}S'_{5}+A_{6}C'_{6}+B_{6}S'_{6}+(A_{15}C'_{15}+\beta B_{15}S'_{15})=0.
\end{cases}
\end{equation}
 Then we assume that the potential $q_i's$ are identical.
 In this way, the characteristic equation $\Phi(\rho)$ is the determinant of a $36\times 36$
matrix.
 The characteristics function, after evaluation of the determinant, is given by
\begin{align*}
\Phi(\rho)=&-[(\alpha^4\be^2+b^2)+(\al^4-\be^4)+(\al^2\be^4+\al^2)]S^6-[(\al^4\be+\be^3)+(\al^3\be^3+\al\be)+(\al^3+\al\be^4)]S^6(2\\
&-36CS')-[(\al^3\be^2+\al\be^2)+(\al^3\be+\al\be^3)+(\al^2\be^3+\al^2\be)](-4374 C^3 (S')^3+2430 C^2 (S')^2\\
&-270 C S'+8)-\al^2\be^2 (531441 C^6 (S')^6-1062882 C^5 (S')^5+728271 C^4 (S')^4-204120 C^3 (S')^3\\
&+21627 C^2 (S')^2-918 C S'+15)
\end{align*}
After simiplification,
\begin{align*}
\Phi(\rho)
=&-\al^2\be^2S^6\{[(\al^2+\al^{-2})+(\al^2\be^{-2}+\al^{-2}\be^2)+(\be^2+\be^{-2})]-[\al^2\be^{-1}+\al^{-2}\be)+(\al \be+\al^{-1}\be^{-1})\\
&+(\al\be^{-2}+\al^{-1}\be^2)](2-36CS')\\
&+[(\al+\al^{-1})+(\al\be^{-1}+\al^{-1}\be)+(\be+\be^{-1})](-4374 C^3 (S')^3+2430 C^2 (S')^2-270 C S'+8)\\
&+531441 C^6 (S')^6-1062882 C^5 (S')^5+728271 C^4 (S')^4-204120 C^3 (S')^3\\
&+21627 C^2 (S')^2-918 C S'+15\}\\
=&-\al^2\be^2S^6\{[2\cos2\theta_1+2\cos2(\theta_1-\theta_2)+(2\cos2\theta_2)]+[2\cos(2\theta_1-\theta_2)+2\cos(\theta_1+\theta_2)\\
&+(2\cos(\theta_1-2\theta_2))](2-36CS')\\
&+[2\cos\theta_1+2\cos(\theta_1-\theta_2)+2\cos\theta_2](-4374 C^3 (S')^3+2430 C^2 (S')^2-270 C S'+8)\\
&+531441 C^6 (S')^6-1062882 C^5 (S')^5+728271 C^4 (S')^4-204120 C^3 (S')^3\\
&+21627 C^2 (S')^2-918 C S'+15\}\ .
\end{align*}
\indent
 In case $q$  is even then $C=S'$, and Theorem~\ref{th1.2}(e) holds.

\section{More analysis on the spectrum}
In this section, we study the spectra of these periodic quantum graphs in more detail. That is, we try to understand the behavior of the point spectrum $\sigma_p$ and the absolutely continuous spectrum $\sigma_{ac}$.
  \newtheorem{th7.1}{Theorem}[section]
 \begin{th7.1}
 \label{th7.1}
Assuming all $q_i$'s are identical and even we have the following
\begin{enumerate}
\item[a.] Let $\sigma(H_{trH})$ be the spectrum of the periodic quantum graph associated with truncated hexagonal tiling. Then
$$\sigma_p(H_{trH})=\{\rho^2:S(a,\rho)=0 \mbox{ or } S'(a,\rho)=0 \mbox{ or } S'(a,\rho)=-\frac{2}{3}\}.$$
\item[b.] The point spectrum of the periodic quantum graph associated with snub square tiling ($3,4,3,4$), rhombi-trihexagonal tiling ($3,4,6,4$), snub trihexagonal tiling ($3^4,6$), and truncated trihexagonal tiling ($4,6,12$) coincide and is given by the set
     $$\{\rho^2:S(a,\rho)=0\}.$$
\end{enumerate}
\end{th7.1}
The above theorem follows immediately from the fact that point spectrum contains exactly those $\lambda=\rho^2$ which do not vary with $(\th_1,\th_2)$.
  \begin{figure}[h!]
\begin{minipage}{.5\linewidth}
\centering
\includegraphics[width=4cm,height=4cm]{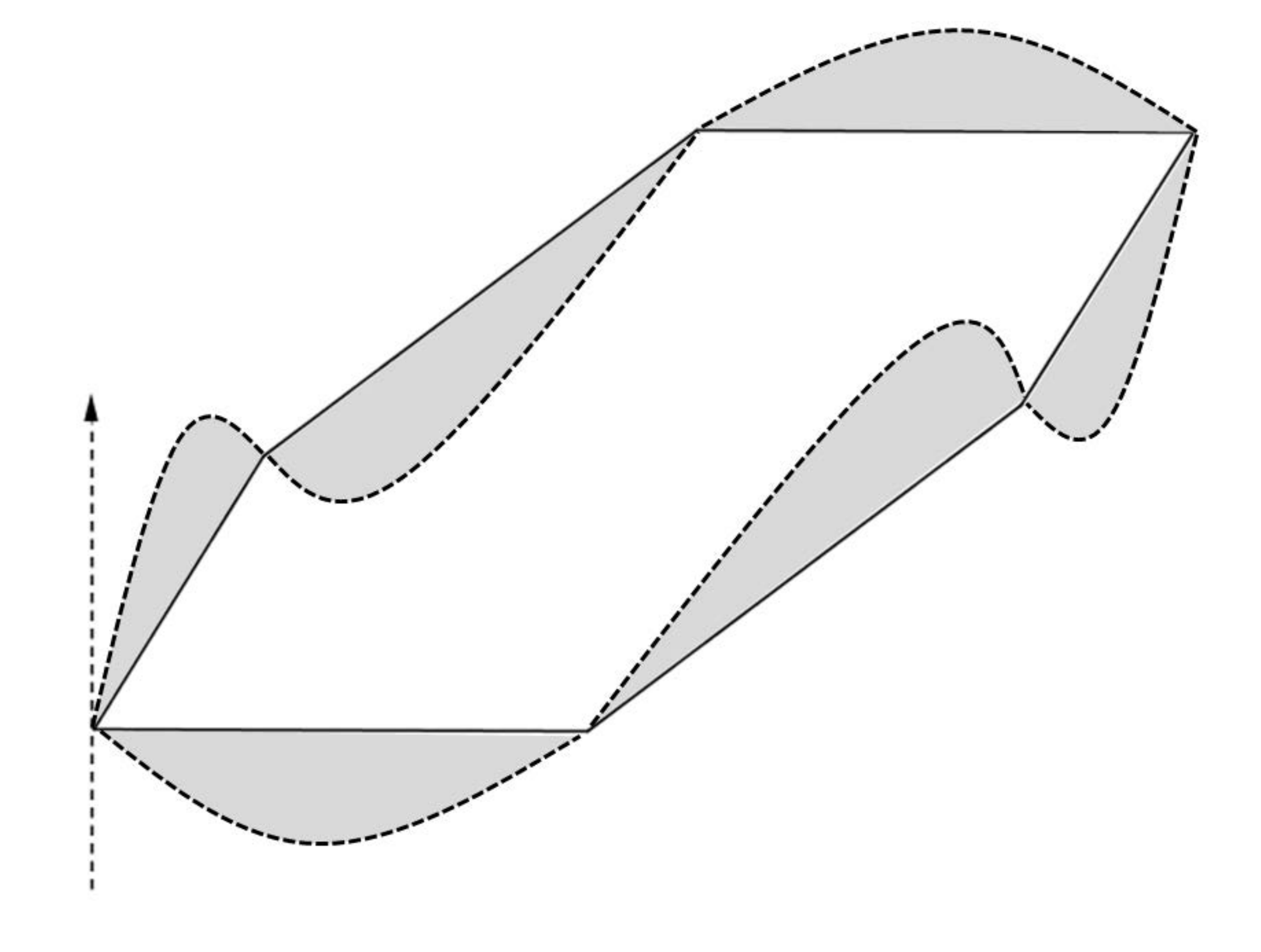}
\end{minipage}%
\begin{minipage}{.5\linewidth}
\centering
\includegraphics[width=4cm,height=4cm]{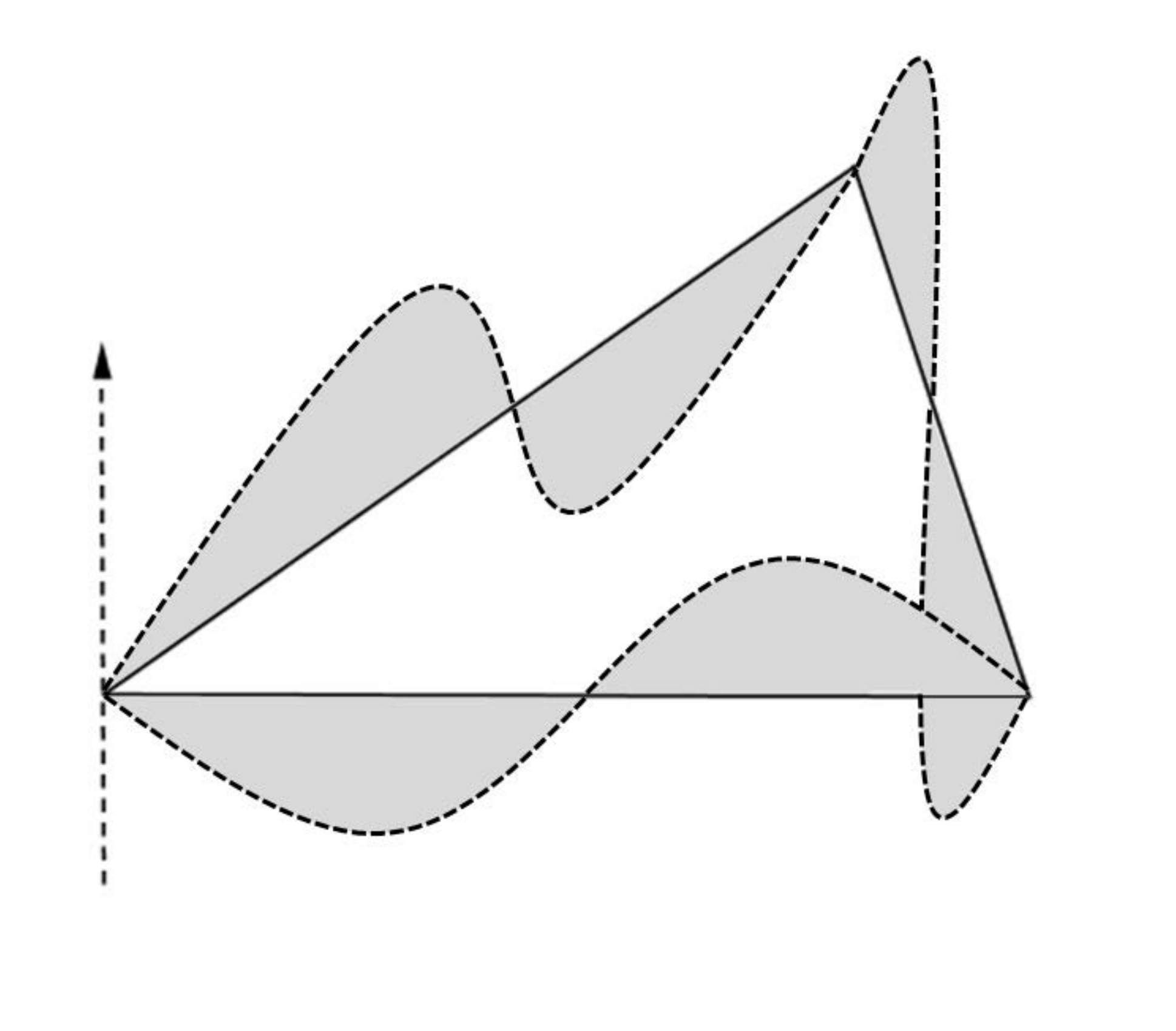}
\end{minipage}\par\medskip
\caption{Eigenfunctions corresponding to the factor $S=0$}
 \label{fig7.1}
 \end{figure}

\begin{figure}[h!]
\begin{minipage}{.5\linewidth}
\centering
\subfloat[]{\label{main:a}\includegraphics[width=6cm,height=4cm]{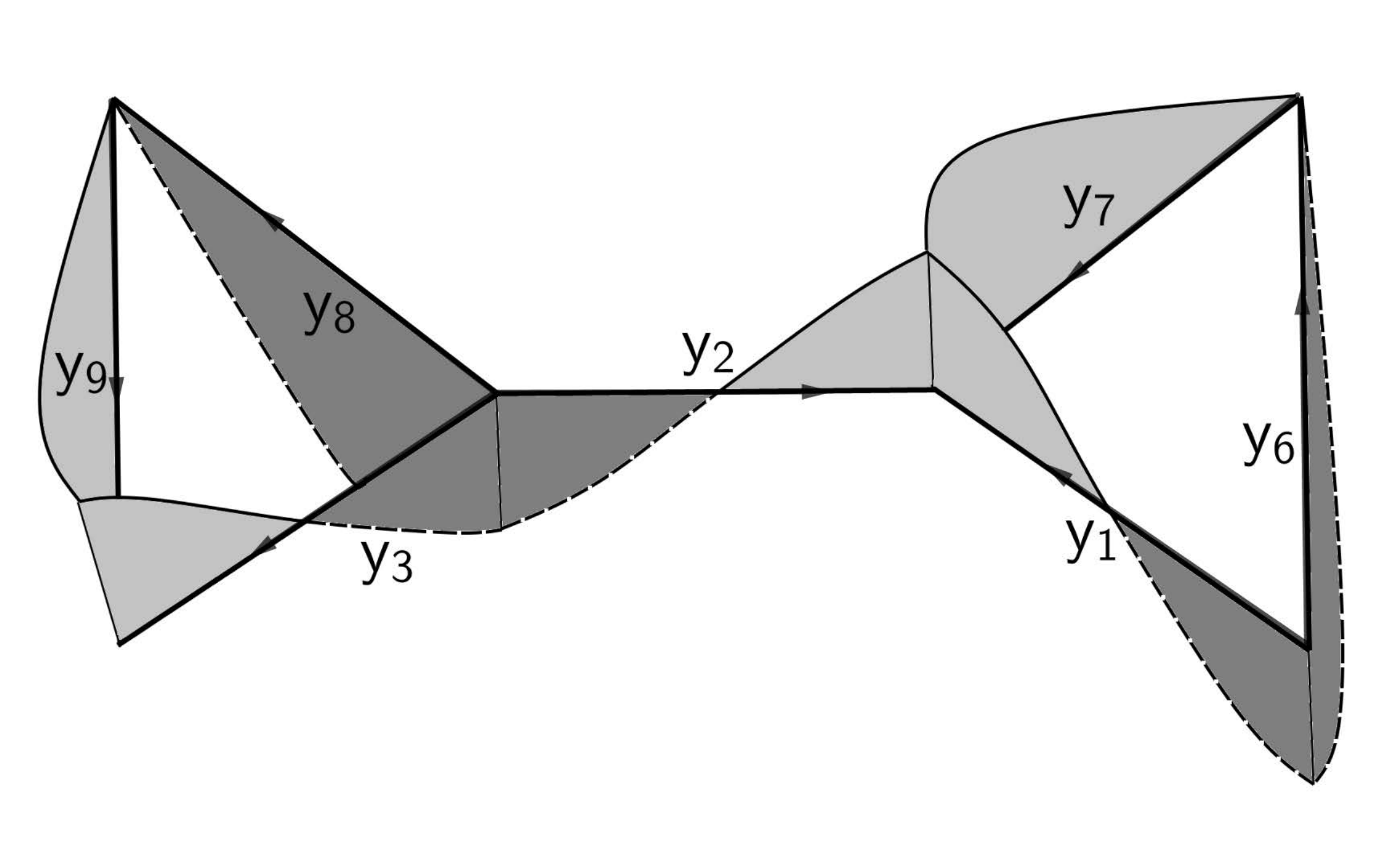}}
\end{minipage}%
\begin{minipage}{.5\linewidth}
\centering
\subfloat[]{\label{main:b}\includegraphics[width=6cm,height=4cm]{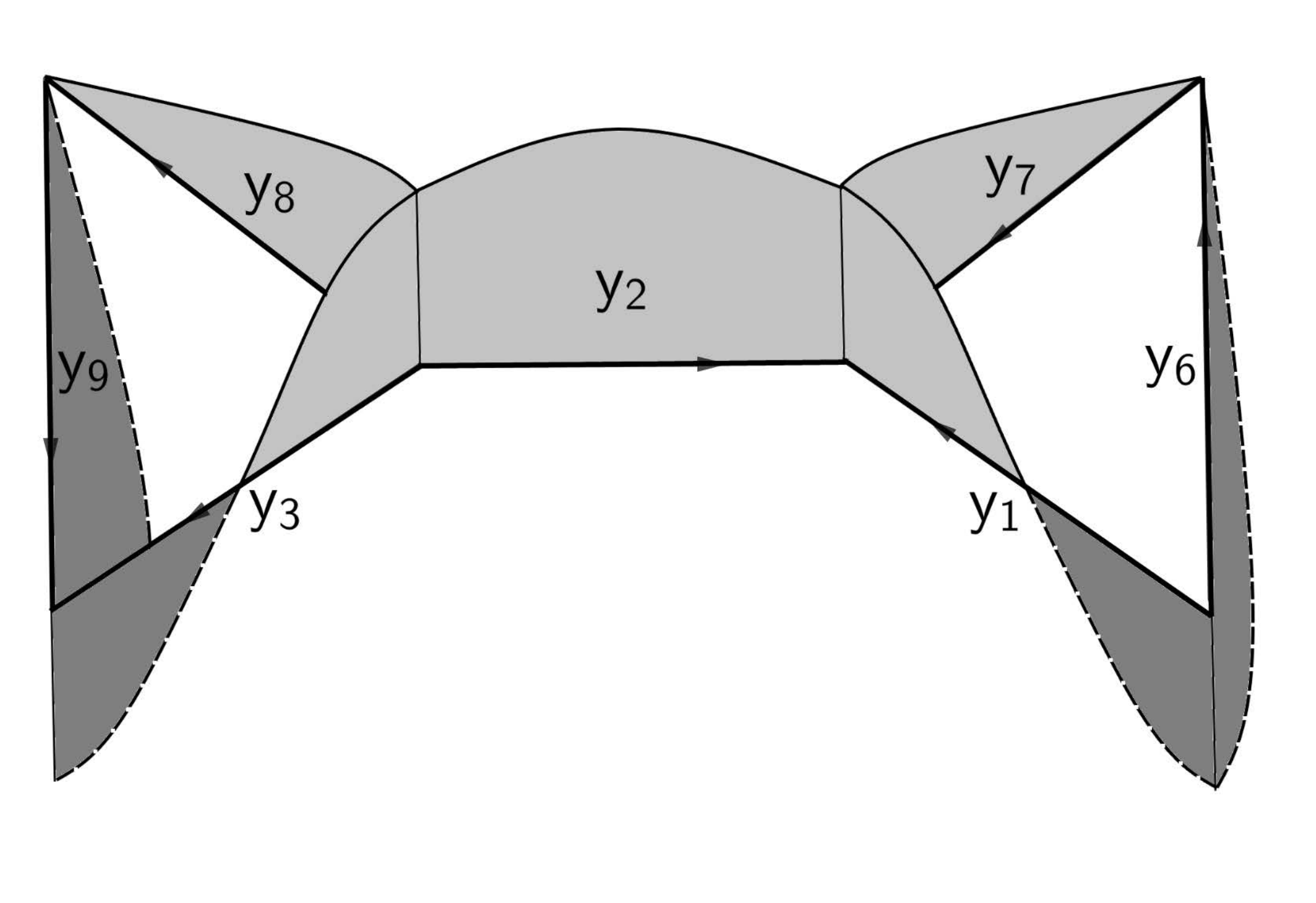}}
\end{minipage}\par\medskip
\caption{Eigenfunctions corresponding to the factors $S'=-2/3$ and $S'=0$}
 \label{fig7.2}
 \end{figure}
 \indent
 Clearly the Bloch variety for the factor $S(a,\rho)=0$ generates the Dirichlet-Dirichlet eigenfunctions along any polygon (triangle, rectangle,
 hexagonal, or dodecagon), and extended by zero to the whole graph, as shown in Fig.\ref{fig7.1}. Truncated hexagonal tiling is interesting in that its
 dispersion relation has two more factors
 $S'(3S'+2)$. The Bloch variety associated with the factor $S'(a,\rho)=-2/3$ generates an eigenfunction with support on the six triangle around the edges of a dodecagon,
 and be extended to the whole graph by zero (cf.\ Fig.\ref{fig2.1} and Fig.\ref{main:a}). In particular, let $f(x)=S(x,\rho)$ such that $\d f'(a)=S'(a,\rho)=\frac{-2}{3}$. Then
 $$
 y_7(x)=f(x)=y_9(x),\quad y_6(x)=-f(a-x)=y_8(x),\quad -y_1(x)=f(x)-f(a-x)=y_3(x)=y_2(x).
 $$
 The pattern goes round the dodecagon and it is easy to see that the Neumann conditions are satisfied at the vertices.

 For the case $S'(a,\rho)=0$, we let $g(x)=S(x,\rho) $ so that $g'(a)=S'(a,\rho)=0$. Define
 $$
 y_7(x)=g(x)=-y_9(x),\qquad y_8(x)=g(a-x)=-y_6(x),\qquad y_1(x)=g(x)-g(a-x)=-y_3(x).
 $$
 It is easy to that $y_2(x)=g(x)+g(a-x)$ is a solution that satisfies $y_2(0)=y_2(a)=g(a)$, so the the Neumann conditions are met
 (cf.\ Fig.\ref{fig2.1} and Fig.\ref{main:b}).

 Hence every element of point spectra has infinitely many eigenfunctions. We conjecture that the functions discussed above are all the eigenfunctions,
 but we cannot prove it.  For the absolutely continuous spectrum for each Archimedean tiling, we have the following theorem.
 \newtheorem{th7.2}[th7.1]{Theorem}
 \begin{th7.2}
 \label{th7.2}
Assuming all $q_i$'s are identical and even, we have the following
\begin{enumerate}
\item[a.] $\sigma_{ac}(H_{trH})=\{\rho^2:S'(a,\rho)\in[-\dfrac{2}{3},0]\cup[\dfrac{1}{3},1]\}$.
\item[b.] $\sigma_{ac}(H_{SS})=\{\rho^2:S'(a,\rho)\in[-\dfrac{3}{5},1]\}$.
\item[c.] $\sigma_{ac}(H_{RTH})=\{\rho^2:S'(a,\rho)\in[-\dfrac{3}{4},1]\}$.
\item[d.] $\sigma_{ac}(H_{STH})=\{\rho^2:S'(a,\rho)\in[\dfrac{-1-\sqrt{3}}{5},1]\}$.
    \item[e] $\sigma_{ac}(H_{trTH})=\{\rho^2:S'(a,\rho)\in[-1,-\dfrac{\sqrt{3}}{3}]\cup[-\dfrac{1}{3},\dfrac{1}{3}]\cup[\dfrac{\sqrt{3}}{3},1]\}$.
\end{enumerate}
\end{th7.2}
\begin{proof}
\begin{enumerate}
\item[(a)] Let $x=S'(a,\rho)$. From Theorem \ref{th1.2}(a), $\rho^2\in \sigma_{ac}(H_{trH})$ if and only if
 $$
 81x^4-54x^3-45x^2+18x-8\cos(\dfrac{\th_1}{2})\cos(\dfrac{\th_2}{2})\cos(\dfrac{\th_1-\th_2}{2})+8=0,
 $$
 or equivalently,
 \begin{equation}
 L(x):= 81x^4-54x^3-45x^2+18x= R(x):= 8\omega-8,
 \label{eq7.1}
 \end{equation}
 where $\d \omega=\cos(\dfrac{\th_1}{2})\cos(\dfrac{\th_2}{2})\cos(\dfrac{\th_1-\th_2}{2})$. Recall the trigonometric identity (see \cite{KP} and
 \cite[Eq.(1.6)]{LJL})
 \begin{eqnarray}
 1+8\omega &=& 3+2(\cos \th_1+\cos\th_2+\cos(\th_1-\th_2)\label{eq7.3}\\
  &=&\left| 1+\rme^{i \th_1}+\rme^{i\th_2}\right|^2 \in [0,9],\nonumber
  \end{eqnarray}
 which implies that $R(x)\in [-9,0]$ for any $x$. Now the function $L$ is a quartic polynomial with with critical points at
 $\d \frac{1}{6},\ \frac{1\pm \sqrt{13}}{6}$. Its
 global minimum points are located at $\d (\dfrac{1-\sqrt{13}}{6},-9)$ and $\d(\dfrac{1+\sqrt{13}}{6},-9)$. It also has a local maximum point
 at $\d (\frac{1}{6},\frac{25}{16})$, and approaches to $\infty$ as $|x|\to\infty$. Also $L(x)>0$ on the intervals $\d (-\infty,-\frac{2}{3}),\
 (0,\frac{1}{3})$, and $(1,\infty)$, while $R$ is a constant function ranging from $-9$ to $0$. It is thus obvious that the dispersion relation
 \eqref{eq7.1} is satisfied iff $\d x=S'(a,\rho)\in [-\frac{2}{3},0] \cup [\frac{1}{3},1]$, whence $\rho^2\in \sigma_{ac}(H_{trH})$.

  \item[(b)] Let $x=S'(a,\rho)$, and $c=\cos\th_1$, $d=\cos\th_2$. From \ref{th1.2}(b), we have
  \begin{equation}
  625 x^4-250 x^2-40 x+1-100(c+d) x^2-40(c+d+cd)x-4(c+d+4cd-c^2-d^2)=0.
  \label{eq7.2}
  \end{equation}
  If $d=1$, this is equivalent to $L(x)=R_{c}(x)$, with
  $$
  L(x):=625 x^4-250 x^2-40 x+1;\qquad R_{c}(x):= 100(c+1)x^2+40(2c+1)x+4(5c-c^2).
  $$
  It is easy to see that $L$ is a quartic polynomial with two global minimum points at $\d x=-\frac{2}{5},\frac{1+\sqrt{2}}{5}$, one local maximum point at $\d
  \frac{1-\sqrt{2}}{5}$, and approaches to $\infty$ as $|x|\to\infty$. Also
  $$
  R_{-1}(x)=-40 x-24,\qquad R_1(x)=200 x^2+120 x+16.
  $$
  It is easy to show that the graph of $R_{-1}$ lies below that of $L$ and is tangent at the points $\d x=\pm \frac{1}{\sqrt{5}}$. On the other hand,
  the graph of $R_1$ intersects that of $L$ at the points $\d (-\frac{3}{5}, 16)$ and $(1,336)$, and is above the graph of $L$ on this interval $\d [-\frac{3}{5},1]$.
  Since $R_{c}$ is continuous in $c$, we conclude that if $\d x=S'(a,\rho)\in [-\frac{3}{5},1]$, then $\rho^2\in \sigma_{ac}(H_{SS})$.
 Now, let $S'(a,\rho)=1+t\ (t>0)$. Then \eqref{eq7.2} yields that for all $c,d\in [-1,1]$,
 \begin{eqnarray*}
 \lefteqn{625t^4+2500t^4+100(35-c-d)t^2+40(49-6c-6d-cd)t}\\
  &&+4(86-36c-36d+c^2+d^2-14cd)\ >\ 0\hskip1in,
  \end{eqnarray*}
  because the constant term is $4(84-36(c+d)+(c-d)^2-12cd\geq 0$.
 Also, if $S'(a,\rho)=-\dfrac{3}{5}-t$  for $t>0$ then \eqref{eq7.2} yields that for all $c,d\in [-1,1]$,
 \begin{eqnarray*}
 \lefteqn{ 625t^4+1500t^3+(1100-100(c+d))t^2+(280-80(c+d)+40cd)t}\\
  &&+16-16(c+d)+4c^2+4d^2+8cd \ >\  0\hskip2in,
  \end{eqnarray*}
  because the constant term is $4(c+d-2)^2\geq 0$.
 Therefore, $\rho^2\in \sigma_{ac}(H_{SS})$ if and only if $S'(a,\rho)\in[-\dfrac{3}{5},1]$.

  \item[(c)] As before, let $x=S'(a,\rho)$. From Theorem \ref{th1.2}(c), the absolutely continuous spectrum $\sigma_{ac}(H_{RTH})$ is characterized by
  \begin{eqnarray}
  0 &=& 2048 x^6-1536 x^4-128 x^3+192 x^2-3-2(64 x^3+32 x^2-1)(\cos\th_1+\cos\th_2+\cos(\th_1+\th_2))\nonumber\\
  && +(\cos 2\th_1+\cos 2\th_2+\cos2(\th_1+\th_2))-2(\cos(\th_1+2\th_2)+\cos(2\th_1+\th_2)+\cos(\th_2-\th_1)\nonumber\\
  &=& 4(512 x^6-384x^4-128 \omega_2 x^3+64(1-\omega_2)x^2+2\omega_2-2\omega_3+\omega_1-1),\label{eq7.5}
  \end{eqnarray}
  with
  \begin{center}
  $\omega_1=\cos\theta_1\cos\theta_2 \cos(\theta_1+\theta_2);\quad \omega_2=\cos(\dfrac{\theta_1}{2})\cos(\dfrac{\theta_2}{2}) \cos(\dfrac{\theta_1+\theta_2}{2});$\\
  $\omega_3:=\cos(\dfrac{2\th_1+\th_2}{2})\cos(\dfrac{\th_2-\th_1}{2})\cos(\dfrac{\th_1+2\th_2}{2}).$
  \end{center}
 and because by \eqref{eq7.3},
 $$
 \cos(\th_1+2\th_2)+\cos(2\th_1+\th_2)+\cos(\th_2-\th_1)=4\omega_3-1.
 $$
  Let us first look at the special case when $\th_1=\th_2=\th$, and $c:=\cos\th$ and $x:=S'(a,\rho)$. Then
  $$
  \omega_2=\dfrac{c+c^2}{2},\quad \omega_1=2c^4-c^2,\quad \omega_3=\dfrac{1+4c^3-3c}{2}.
  $$
  Thus the above dispersion relation becomes
 \begin{equation}
 L(x):= 512x^6-384x^4+64x^2\quad =\quad 2(c+1)(32 c x^3+16 c x^2)-(c-1)^3):= R_{c}(x) \label{eq7.4}
 \end{equation}
 for $c\in [-1,1]$. Here
 $L$ is an even polynomial with global minimum points at $x=\pm\dfrac{\sqrt{9+3\sqrt{3}}}{6}$,  two local maximum points at $x=\pm\dfrac{\sqrt{9-3\sqrt{3}}}{6}$, and one local minimum point at $x=0$.  Also $L(x)$ approaches to $\infty$ as $|x|\to\infty$. The right hand side $R_{c}$ is in general a cubic polynomial. We shall study three cases:
 $$
 R_{-1}\equiv 0,;\quad R_1(x)=128 x^3+64 x^2;\quad R_{-1/2}(x)=-16x^3-8x^2+\frac{27}{8}.
 $$
 The graph of $R_{-1}$ is a line that intersects that of $L$ at $x=0,\pm\sqrt{\dfrac{1}{2}},\pm\dfrac{1}{2}$. Also $R_{-1}\leq L$ on $\d [0,\frac{1}{2}]$.
 The graph of $R_{-1/2}$ intersects that of $L$ at $x=-\dfrac{3}{4},\dfrac{1}{4}$.  Furthermore $R_{-1/2}\geq L$ on $[-\dfrac{3}{4},\dfrac{1}{4}]$, while
 $R_{-1/2}\leq L$ on $\d [\frac{1}{4},1]$. Thus by continuity of $R_{c}$, if $x$ lies in $\d [0,\frac{1}{4}]$, then \eqref{eq7.4} is valid for some
 $\d c\in [-1,-\frac{1}{2}]$.
 On the other hand, $R_1$ and $L$ have the same values at $x=0,1$; and $R_1\leq L$ on $[-1,0]$ while $R_1\geq L$ on $[0,1]$. Comparing with the properties of $R_{-1/2}$, we
 conclude that for $x\in [-\dfrac{3}{4},0]\cup[\dfrac{1}{4},1]$, there exists $c\in [-1/2,1]$ such that \eqref{eq7.4} is valid.
 Therefore $\d \{\rho^2: S'(a,\rho)\in [-\frac{3}{4},1]\}\subset \sigma_{ac}(H_{RTH})$.

  We shall see that this is indeed the entire range of $S'(a,\rho)$, for any $(\th_1,\th_2)\in [-\pi,\pi]^2$.
 Suppose $S'(a,\rho)=1+t$ for $t>0$ then \eqref{eq7.5} reduces to
\begin{eqnarray*}
  \lefteqn{512t^6+3072t^5+7296t^4+128(68-\omega_2)t^3+64(85-7\omega_2)t^2}\\
  &&+(1664-512d)t+191+\omega_1-190\omega_2-2\omega_3>0,\hskip1in
\end{eqnarray*}
because all coefficients are nonnegative (Please refer to Appendix B for detailed argument).

Also, if $S'(a,\rho)=-\dfrac{3}{4}-t$ for some $ t>0$, then \eqref{eq7.5} implies
\begin{eqnarray*}
 \lefteqn{512t^6+2304t^5+3936t^4+32(99+4\omega_2)t^3+(1198+224\omega_2)t^2}\\
 &&+3(59+40\omega_2)t+\omega_1-2\omega_3+20\omega_2+\frac{37}{8}>0\hskip1.5in.
 \end{eqnarray*}
 because all coefficients are nonnegative (cf.\ Appendix B for detail).
 Therefore, $\rho^2\in \sigma_{ac}(H_{RTH})$ if and only if $S'(a,\rho)\in[-\dfrac{3}{4},1]$.
\item[(d)] Let $\omega_1,\omega_2,\omega_3$ be as defined above.
 By the trigonometric identity \eqref{eq7.3}, the dispersion relation given by Theorem~\ref{th1.2}(d) is simplified to
 \begin{eqnarray*}
 \lefteqn{S^9\, [15625 (S')^6-9375 (S')^4-(4000\omega_2+1000)(S')^3-(2400\omega_2-1275)(S')^2}\\
 && -(240\omega_2+80\omega_3-200)S'+8\omega_1+32\omega_2-32\omega_3-13]=0.\hskip1in
 \end{eqnarray*}
  Let $x=S'(a,\rho)$. The absolutely continuous spectrum $\sigma_{ac}(H_{STH})$ is characterized by
  \begin{eqnarray}
  \lefteqn{15625 x^6-9375 x^4-(4000 {\omega}_2+1000)x^3-(2400\omega_2-1275)x^2}
   \nonumber\\
   &&-(240\omega_2+80\omega_3-200)x+8\omega_1+32\omega_2-32\omega_3-13=0.\hskip1in\label{eq7.6}
  \end{eqnarray}
  Now let $\th_1=\th_2=\th$, $c:=\cos\th$. Then \eqref{eq7.6} becomes
   \begin{eqnarray}
   L(x)&:=&15625x^6-9375x^4\nonumber\\
   &=&(2000c^2 + 2000 c + 1000) x^3 + (1200 c^2 + 1200 c -
    1275) x^2\nonumber \\
    &&+ (160 c^3 + 120 c^2 - 160) x - 16 c^4 + 64 c^3 -8 c^2 - 64 c + 29\label{eq7.7}.
   \end{eqnarray}
    Let $R_{c}$ denote the polynomial on the right hand side. We study $R_1$ and $R_{-1/2}$ in more detail.
Since
  $$
  L(x)-R_{-0.5}(x)=25(25x^2+10x-2)(5x^2-x-1)^2,
   $$
   then $L\leq R_{-0.5}$ on $\d [\dfrac{-1-\sqrt{3}}{5},-\dfrac{1}{5}]$ and $L\geq R_{-0.5}$ on $[\dfrac{-1+\sqrt{3}}{5},1]$. Also,
   $$L(x)-R_1(x)=5(x-1)(5x+1)^5.
   $$
    Hence $L\geq R_1$ on $[\dfrac{-1-\sqrt{3}}{5},-\dfrac{1}{5}]$ and $L\leq R_{-0.5}$ on $[\dfrac{-1+\sqrt{3}}{5},1]$.  Since $R_{c}$ is continuous in $c$,
    it follows that \eqref{eq7.6} is valid for some $c\in[-1,1]$ when $x\in[\dfrac{-1-\sqrt{3}}{5},-\dfrac{1}{5}]\cup[\dfrac{-1+\sqrt{3}}{5},1]$.

     \quad Now, for the interval $[-\dfrac{1}{5},\dfrac{-1+\sqrt{3}}{5}]$, we let $x=-\dfrac{1}{5}+t$ for some $t\in[0,\dfrac{2}{5}]$. Then
\begin{center}
$f(x):=L(x)-R_{c}(x)=15625t^6-18750t^5+(4000-2000 c-2000 c^2)t^3+(120c^2-160 c^3+240 c-200)t+16 c^4-32 c^3+32c-16$
\end{center}
 Let $y=5t$. Then
\begin{eqnarray*}
 \lefteqn{g(y,c):=f(x)=y^6-6y^5+(32-16 c-16 c^2)y^3+(24c^2-32 c^3+48c-40)y}\hskip1in\\
 &&+16 c^4-32 c^3+32 c-16\hskip1.5in
\end{eqnarray*}

  To show that \eqref{eq7.6} is valid for some $c\in[-1,1]$ when $x\in[-\dfrac{1}{5},\dfrac{-1+\sqrt{3}}{5}]$, since $L(x)-R_1(x)\le0$ for all $x\in[-\dfrac{1}{5},\dfrac{-1+\sqrt{3}}{5}]$, it suffices to show that for any $y\in[0,2]$ we have $g(y,c)\ge0$ for some $c\in[-1,1]$. Now, for any $y_0\in[0,2]$ consider $c_0=\dfrac{3y_0+4-y_0\sqrt{8y_0+9}}{4}$. One can easily check that $c_0\in[0,1]$. Substituting $(y_0,c_0)$ we have
 $$
 g(y_0,c_0)=-\dfrac{1}{2}y_0^3(y_0+2)(6y_0^2+36y_0+27-(8y_0+9)\sqrt{8y_0+9})\ge0
 $$
 The above inequality is due to the fact that $h(y):=6y^2+36y+27-(8y+9)\sqrt{8y+9}$ has 0 as its maximum value in $[0,2]$.\\

 Now, suppose $S'(a,\rho)=1+t$ ($t>0$), then the left hand side of \eqref{eq7.6} yields
 \begin{eqnarray*}
  \lefteqn{15625t^6+93750t^5+225000t^4+(274000-4000\omega_2)t^3+(176400-14400\omega_2)t^2}\\
  && +(56000-17040\omega_2-80\omega_3)t+8(839+\omega_1-14\omega_3-826\omega_2)>0
 \end{eqnarray*}
 Also, if $S'(a,\rho)=\dfrac{-1-\sqrt{3}}{5}-t$ for some $ t>0$, then
 \begin{eqnarray*}
 \lefteqn{15625t^6+18750(\sqrt{3}+1)t^5+9375(2\sqrt{3}+3) t^4+500(8\omega_2+15\sqrt{3}+37)t^3}\\
 &&+75(45+32\sqrt{3}\omega_2+28\sqrt{3}) t^2+(80\omega_3+1200 \omega_2+270\sqrt{3}+70)t\hskip1in\\
 &&+(8\omega_1+16(\sqrt{3}-1)\omega_3+16(1+3\sqrt{3})\omega_2-10\sqrt{3}+19)>0.
\end{eqnarray*}
Therefore, $\d \sigma_{ac}(H_{STH})=\{\rho^2:S'(a,\rho)\in[\dfrac{-1-\sqrt{3}}{5},1]\}$.

\item[(e)]
Let $x=:9(S'(a,\rho))^2$, and
 \begin{center}
 $\widetilde{\omega}_1:=\cos\th_1\cos\th_2\cos(\th_1-\th_2),\qquad \widetilde{\omega}_2:=\cos\dfrac{\theta_1}{2}\cos\dfrac{\theta_2}{2}\cos(\dfrac{\theta_1-\theta_2}{2}),$\\
 $\widetilde{\omega}_3:=\cos(\dfrac{2\th_1-\th_2}{2})\cos(\dfrac{\th_1+\th_2}{2})\cos(\dfrac{\th_1-2\th_2}{2}).$
  \end{center}
  Then from Theorem \ref{th1.2}(e) and the trigonometric identity \eqref{eq7.3}, the absolutely continuous spectrum $\sigma_{ac}(H_{trTH})$ is characterized by
 \begin{eqnarray}
 \lefteqn{x^6-18x^5+111x^4-(48\widetilde{\omega}_2+268)x^3+(240\widetilde{\omega}_2+207)x^2}\nonumber\\
 &&-(32\widetilde{\omega}_3+240\widetilde{\omega}_2+34)x+8\widetilde{\omega}_1+64\widetilde{\omega}_2+16\widetilde{\omega}_3-7=0.\label{eq7.8}
 \end{eqnarray}
 We want to show that this is valid if and only if $\d S'(a,\rho)\in [-1,\frac{-1}{\sqrt{3}}]\cup [\frac{-1}{3},\frac{1}{3}]\cup [\frac{1}{\sqrt{3}},1]$,
 or equivalently $x\in [0,1]\cup [3,9]$.
 Now let $\th=\th_1=-\th_2$, $c:=\cos\th$. Then
 $$
 \widetilde{\omega}_1=2 c^4-c^2;\quad\widetilde{\omega}_2=\frac{1}{2}(c^2+c);\quad \widetilde{\omega}_3=\frac{1}{2}(4 c^3-3 c+1),
 $$
 and \eqref{eq7.8} becomes $L(x)=R_{c}(x)$, where
 \begin{eqnarray*}
 L(x)&:=& x^6-18x^5+111x^4\\
 R_{c}(x)&:=& (24 c^2+24 c+268)x^3-(120 c^2+120 c+207)x^2+(64 c^3+120 c^2+72 c+50)x\\
  &&-16 c^4-32 c^3-24 c^2-8 c-1.
  \end{eqnarray*}
  We study two special $R_{c}$'s, namely
  $$
   R_1(x)=316 x^3-447x^2+316 x-81, \qquad R_{-1/2}(x)=262 x^3-177 x^2+36 x.
   $$
   It is routine to see that
   \begin{eqnarray*}
   L(x)-R_1(x)&=&(x-1)^3(x-3)^2(x-9)\\
   L(x)-R_{-1/2}(x)&=&(x^2-7x+3)^2(x-4)x.
   \end{eqnarray*}
  Hence $R_{-1/2}(x)\leq L(x)\leq R_1(x)$ for all $x\in [0,1]\cup [4,9]$. By continuity,
  $[0,1]\cup [4,9]$ is a subset of the characterization domain.

  Next, we show that $[3,4]$ is also a subset. Let $x=3+t, t\in[0,1]$. Then we have
 \begin{eqnarray*}
f(t,c)&:=&L(x)-R_{c}(x)\\
  &=&t^6-24t^4-(2 c^2+24c+16)t^3+(144-96 c-96c^2)t^2\\
  &&+(112-48 c^2-64 c^3)t+16 c^4-160 c^3+96c^2+224c-176.
\end{eqnarray*}
Since $L(x)-R_1(x)\le0$,  for all $x\in[3,4]$, it suffices to show that for any $t\in[0,1]$ we have $f(t,c)\ge0$
for some  $c\in[-1,1]$. Now, for $t_0\in[0,1]$, let $c_0=\dfrac{3t_0+8-(t_0+2\sqrt{3t_0+9})}{2}$. It is easy to check that $c_0\in[\dfrac{11-6\sqrt{3}}{2},1]$. Substituting this $(t_0,c_0)$ we get
$$
f(t_0,c_0)=-4(t_0+2)^3(t_0+3)(2t_0^2+27t_0+54-2(3t_0+9)^{3/2})\ge0.
$$
The above inequality is due to the fact that $g(t)=2t^2+27t+54-2(3t+9)^{3/2}$ has 0 as its maximum value in $t\in[0,1]$.\\
 \indent
 Now let $x=-t, t>0$, then \eqref{eq7.8} will be
 $$
t^6+18t^5+111t^4+(48\widetilde{\omega}_2+268)t^3+(240\widetilde{\omega}_2+207)t^2+(32\widetilde{\omega}_3+240\widetilde{\omega}_2+34)t
+(16\widetilde{\omega}_3+8\widetilde{\omega}_1+64\widetilde{\omega}_2-7)>0,
 $$
 (cf.\ Appendix B). If $x=9+t$ for some $t>0$, then \eqref{eq7.8} becomes
  \begin{center}
  $t^6+36t^5+516t^4+(3728-48\womega_2)t^3+(14112-1056\womega_2)t^2+(26048-7584\womega_2-32\womega_3)t+(8\womega_1-272\womega_3-17648\womega_2+17912)>0$,
  \end{center}
  (cf.\ Appendix B). Also, if $x=2+t$ for $-1<t<1$, then we have
 $$
 t^6-6t^5-9t^4+(60-48\womega_2)t^3+(63-48\womega_2)t^2+(144\womega_2-32\womega_3-118)t+(8\womega_1-48\womega_3+160\womega_2-127)<0,
 $$
 (cf.\ Appendix B). Therefore,  we conclude that $\rho^2\in \sigma_{ac}(H_{trTH})$ if and only if $x\in [0,1]\cup [3,9]$.
\end{enumerate}
\end{proof}
\indent
  In the case $q=0$, then $S'(a,\rho)=\cos(\rho a)$.  Hence we have the following corollary.
\newtheorem{th7.3}[th7.1]{Corollary}
\begin{th7.3}
\label{th7.3}
  When $q=0$, we have
\begin{enumerate}
\item[(a)] $\sigma_{ac}(H_{TH})=\d\bigcup_{k=0}^\infty\left\{\left[(\dfrac{2(k+1)\pi-\xi_1}{a})^2,(\dfrac{2(k+1)\pi+\xi_1}{a})^2\right]\cup\left[(\dfrac{(4k+1)\pi}{2a})^2,(\dfrac{2k\pi+\xi_2}{a})^2\right]\right.$\\
	$\left.\cup \left[(\dfrac{2(k+1)\pi-\xi_2}{a})^2,(\dfrac{(4k+3)\pi}{2a})^2\right]\right\}\cup\left[0,(\dfrac{\xi_1}{a})^2\right]$,\\
  where $\xi_1=\arccos(\dfrac{1}{3})$, $\xi_2=\arccos(-\dfrac{2}{3})$.
\item[(b)]$\sigma_{ac}(H_{SS})=\d\bigcup_{k=0}^\infty\left\{\left[(\dfrac{2k\pi}{a})^2,(\dfrac{\xi+2k\pi}{a})^2\right]\cup\left[(\dfrac{2(k+1)\pi-\xi}{a})^2,
(\dfrac{2(k+1)\pi}{a})^2\right]\right\}$,
where $\xi=\arccos(-\dfrac{3}{5})$.
\item[(c)]$\sigma_{ac}(H_{RTH})=\d\bigcup_{k=0}^\infty\left\{\left[(\dfrac{2k\pi}{a})^2,(\dfrac{\xi+2k\pi}{a})^2\right]\cup\left[(\dfrac{2(k+1)\pi-\xi}{a})^2,
(\dfrac{2(k+1)\pi}{a})^2\right]\right\}$,
where $\xi=\arccos(-\dfrac{3}{4})$.
\item[(d)]$\sigma_{ac}(H_{STH})=\d\bigcup_{k=0}^\infty\left\{\left[(\dfrac{2k\pi}{a})^2,(\dfrac{\xi+2k\pi}{a})^2\right]\cup\left[(\dfrac{2(k+1)\pi-\xi}{a})^2,
(\dfrac{2(k+1)\pi}{a})^2\right]\right\}$,
where $\xi=\arccos(-\dfrac{1+\sqrt{3}}{5})$.
\item[(e)]$\sigma_{ac}(H_{trTH})=\d\bigcup_{k=0}^\infty\left\{\left[(\dfrac{(k+1)\pi-\xi_1}{a})^2,(\dfrac{(k+1)\pi+\xi_1}{a})^2\right]
	\cup\left[(\dfrac{k\pi+\xi_2}{a})^2,(\dfrac{(k+1)\pi-\xi_2}{a})^2\right]\right\}\cup\left[0,(\dfrac{\xi_1}{a})^2\right]$, where $\xi_1=\arccos(\dfrac{\sqrt{3}}{3}), \xi_2=\arccos(\dfrac{1}{3})$
\end{enumerate}
\end{th7.3}
  \section{Concluding remarks}
  As a summary, we have derived the dispersion relations for the periodic quantum graphs associated with all the 11 Archimedean tilings.  We showed that
  they are all of the form
 $S^i S'^j(2S'-1)^k (3S+2)^l p(S',\th_1,\th_2)=0$ for some $i,j,k,l\geq 0$, and   $p$ is a polynomial of $S'$.  The spectrum $\sigma(H)$ is exactly Bloch variety. Here the point spectrum $\sigma_p(H)$ is defined by $S^i(S')^j(2S'-1)^k (3S+2)^l=0$; while the absolutely continuous spectrum $\sigma_{ac}$ is defined by $p(S',\th_1,\th_2)=0$,
 For the 5 Archimedean tilings discussed in this paper, we can find many  eigenfunctions, each of which has infinite multiplicity.
 Some nontrivial eigenfunctions for trihexagonal tiling and truncated hexagonal tiling are also found. The spectra are all of a band and gap structure.  Also their absolutely continuous spectra satisfy $\sig_{ac}\subset \phi^{-1}( [-1,1])$
 where $\phi(\rho)=S'(a,\rho)$.
 Moreover the absolutely continuous spectra corresponding to each of the square tiling, hexagonal tiling, truncated square tiling and truncated trihexagonal tiling
 satisfies  $\sig_{ac}= \phi^{-1}( [-1,1])$.
 As the function behaves asymptotically like $\cos \rho a$, as $\rho\to\infty$, we know that the spectral gaps tends to zero.

 We also remark that for $\sigma_{ac}$, $S'=1$ is always a solution for $p$, for all the 11 Archimedean tilings.
 We also perform numerical checkup for specific values of $\Th$ to check those
dispersion relations derived using symbolic software.  Thus we are confident that the dispersion relations in this paper are correct.
 Thus our method is an efficient and effective one for any complicated periodic quantum graphs.

Recently, Fefferman and Weinstein \cite{FW12} embarked on a systematic study of the case of hexagonal tiling (also called honeycomb lattice).  Based on a PDE approach,
assuming that the wave functions acts on the whole hexagon, they showed that there exists infinitely many Dirac points in the spectrum.  In our case, Dirac points are the points
$\la=\rho^2$ where $\phi(\rho)=S'(a,\rho)=1$ with $\phi'(\rho)=0$, and the $(\la,\th_1,\th_2)$ relationship is asymptotically a cone (called Dirac cone) \cite{K2013}. As explained in \cite{FW12,K2013}, the existence of Dirac points accounts for the electronic properties of graphene and related crystal lattices.
It is desirable to study when there exist infinitely many Dirac points, not only for hexagonal tiling, but also for the other 10 Archimedean tilings.
We shall pursue on this issue later.

   \section*{Acknowledgements}
   \indent
We thank Min-Jei Huang, Ka-Sing Lau, Vyacheslav Pivovarchik and Ming-Hsiung Tsai for stimulating discussions. We also thank the anonymous referees for helpful comments. 
The authors are partially supported by Ministry of Science and Technology, Taiwan, under contract number MOST105-2115-M-110-004.
  \newpage
 \begin{table}[h!]
 \centering
 \caption{Table of 11 Archimedean tilings}
 $\begin{array}{|c||c|c|c|}
 \hline
 \text{Name} & \text{Triangular (T)} & \text{{\small Snub trihexagonal (ST)}} & \text{{\small Elongated triangular (ET)}} \\ \hline
 \text{Notation} & \left(3^6\right) & \left(3^4,6\right) & \left(3^3,4^2\right) \\ \hline
 \  & \includegraphics[width=3cm]{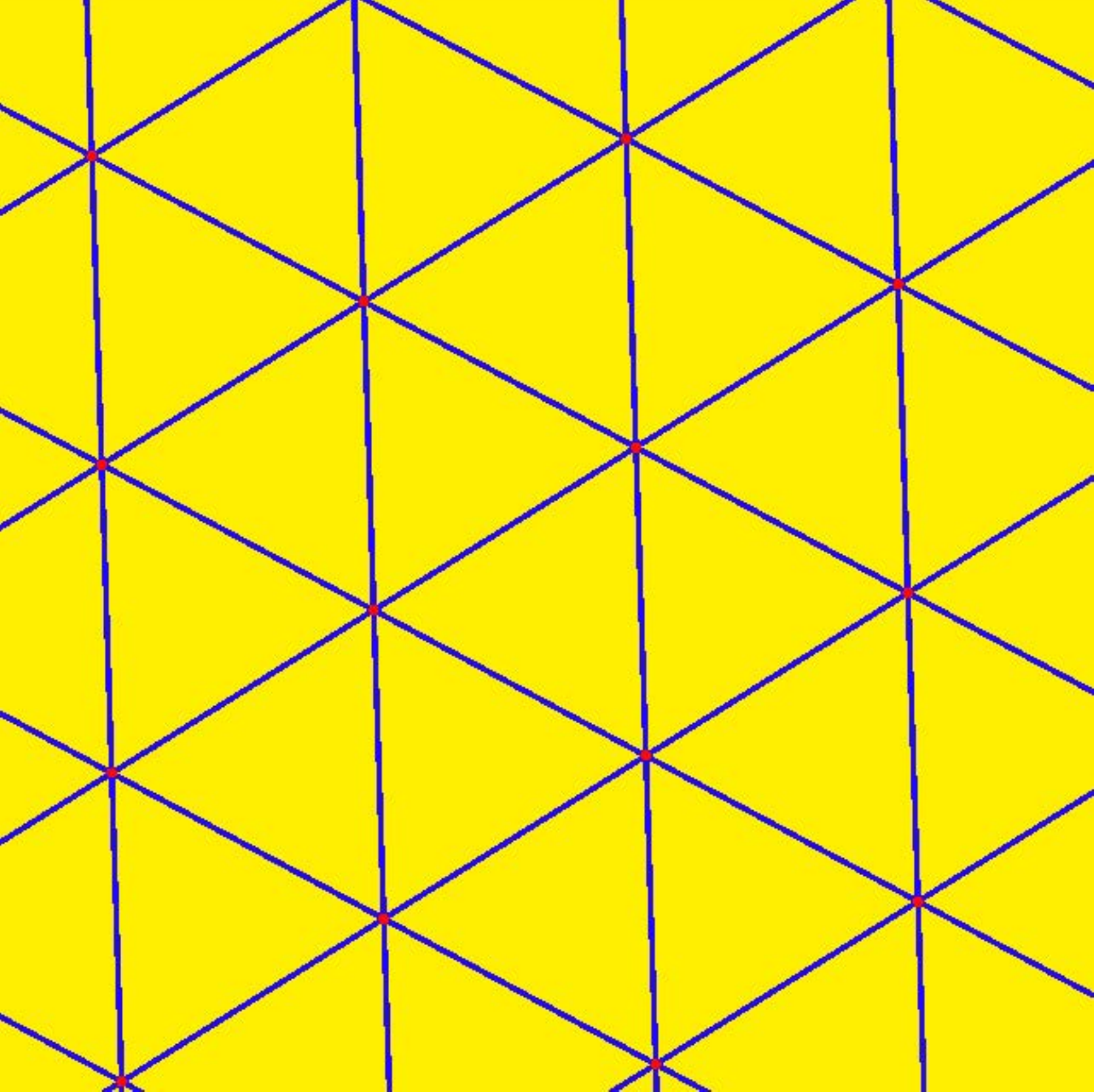} & \includegraphics[width=3cm]{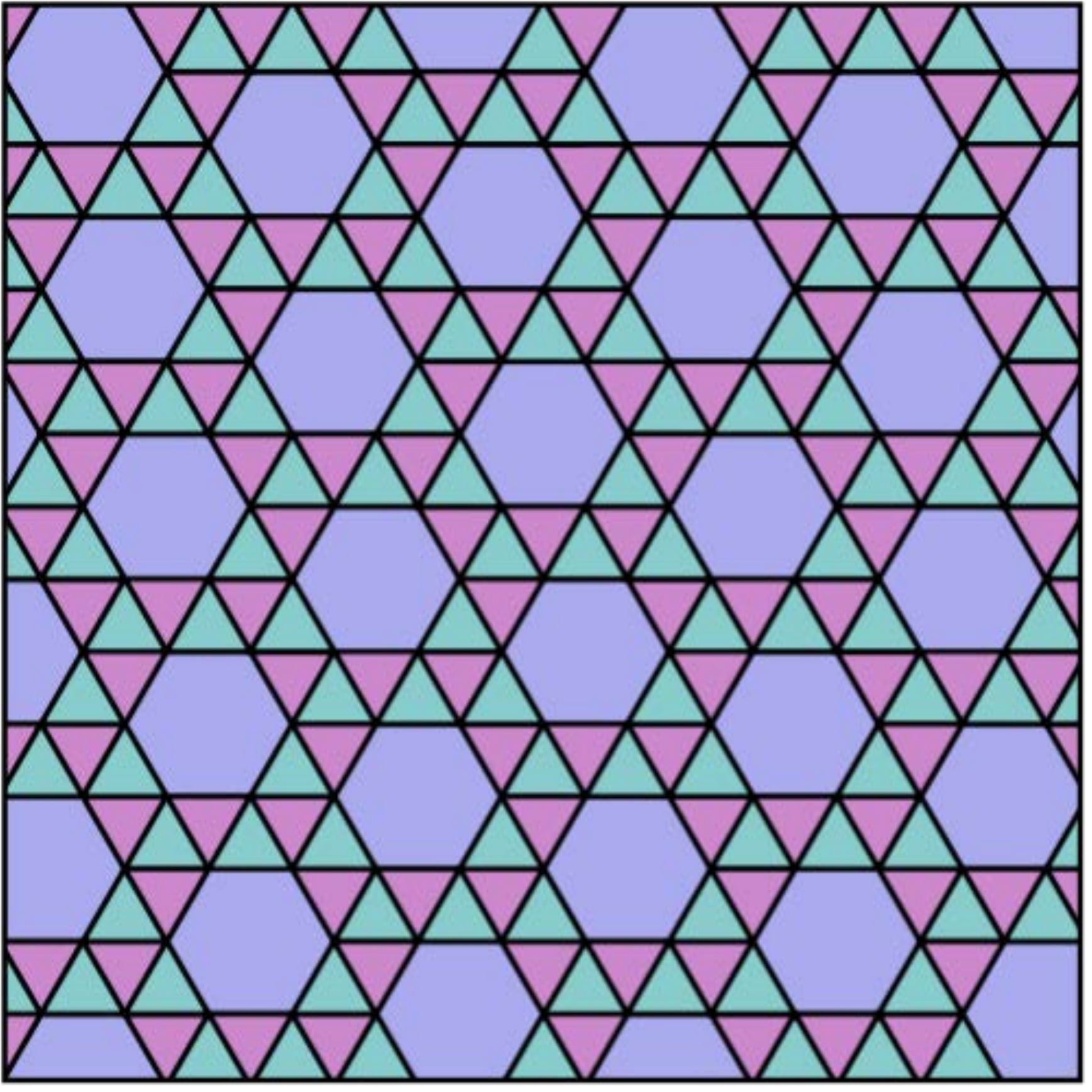} & \includegraphics[width=3cm]{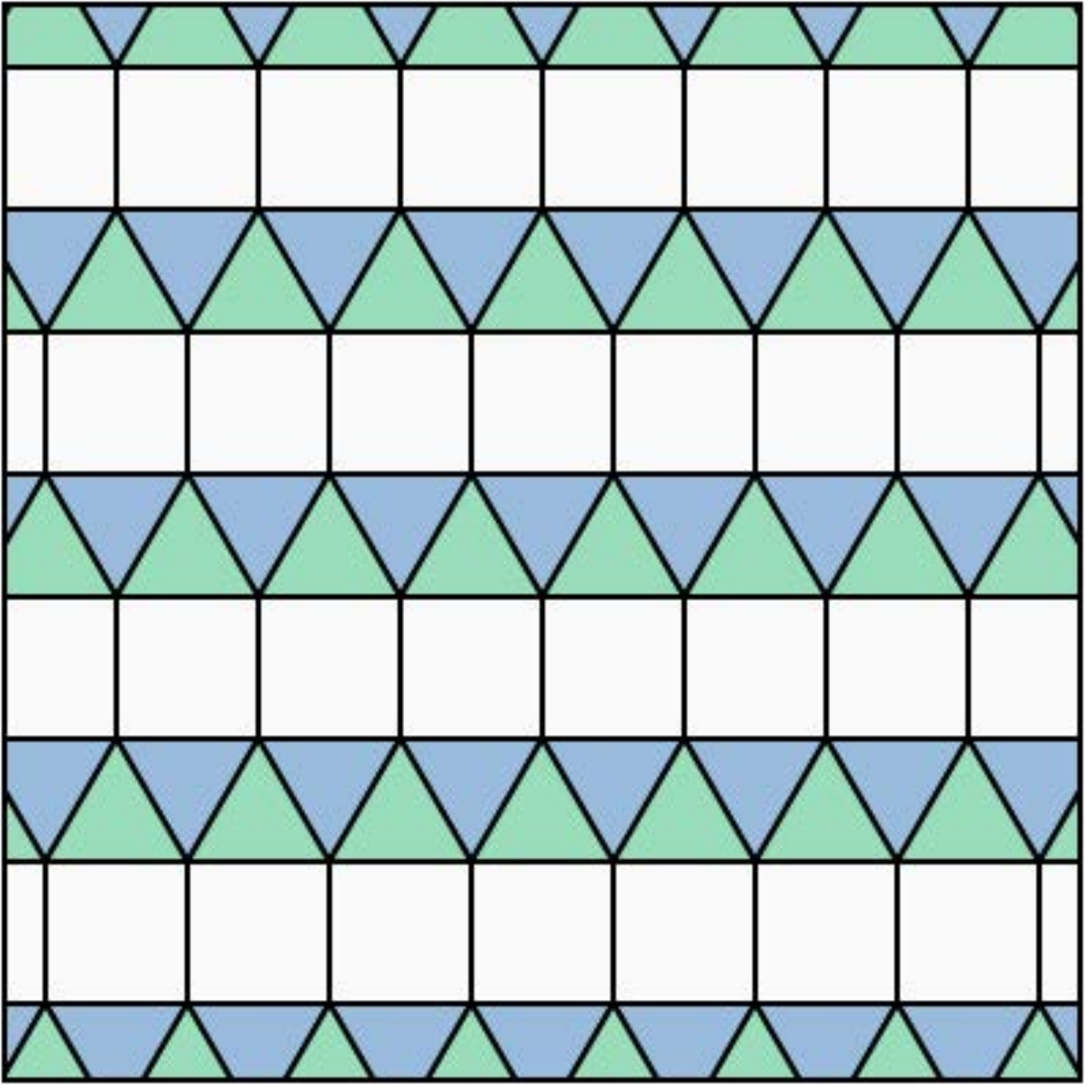}  \\ \hline \hline
 \text{Name} & \text{Snub square (SS)} & \text{Trihexagonal (TH)} & \text{{\small Rhombi-trihexagonal (RTH)}} \\ \hline
 \text{Notation} & \left(3^2,4,3,4\right) & (3, 6, 3, 6) & (3,4,6,4) \\ \hline
 \  & \includegraphics[width=3cm]{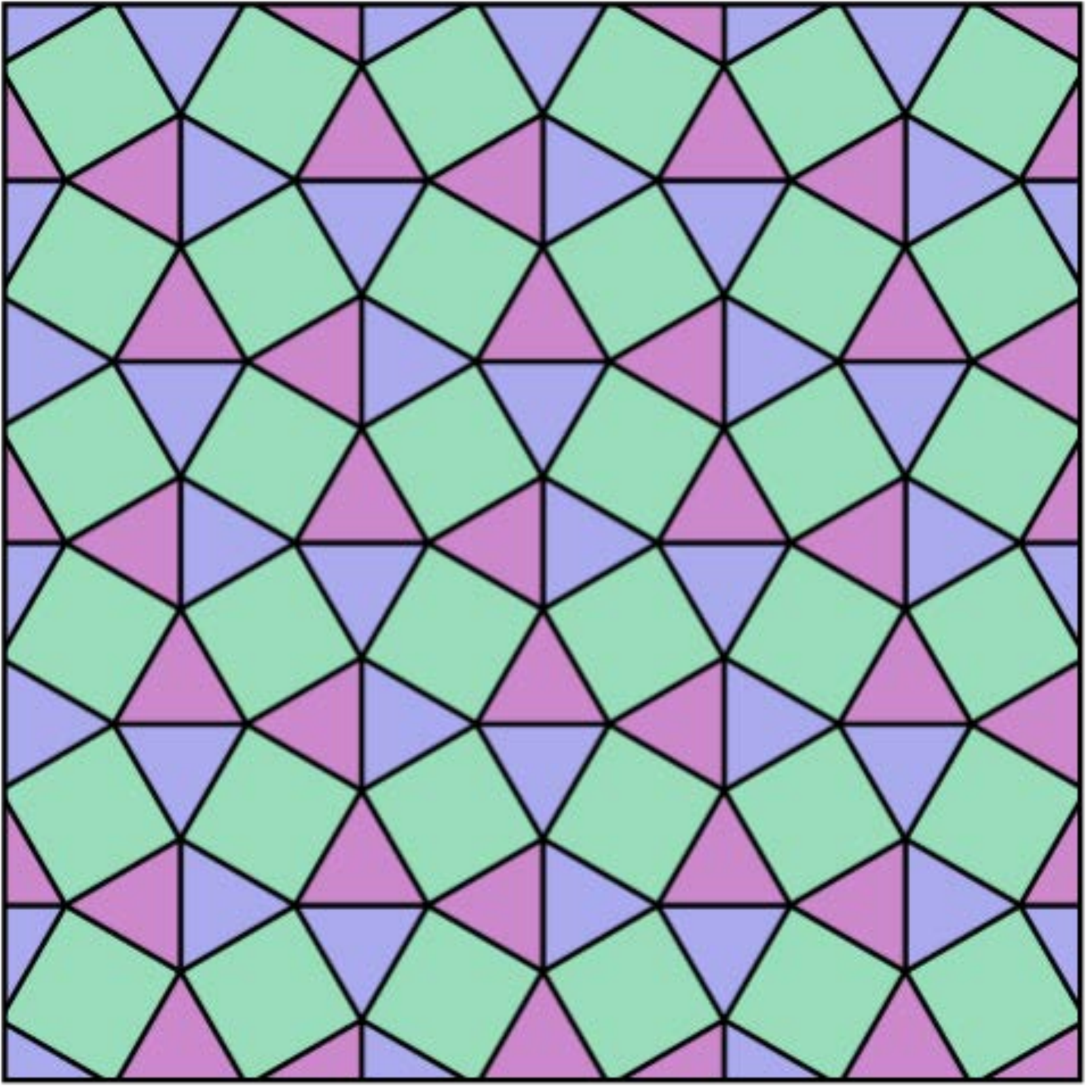} & \includegraphics[width=3cm]{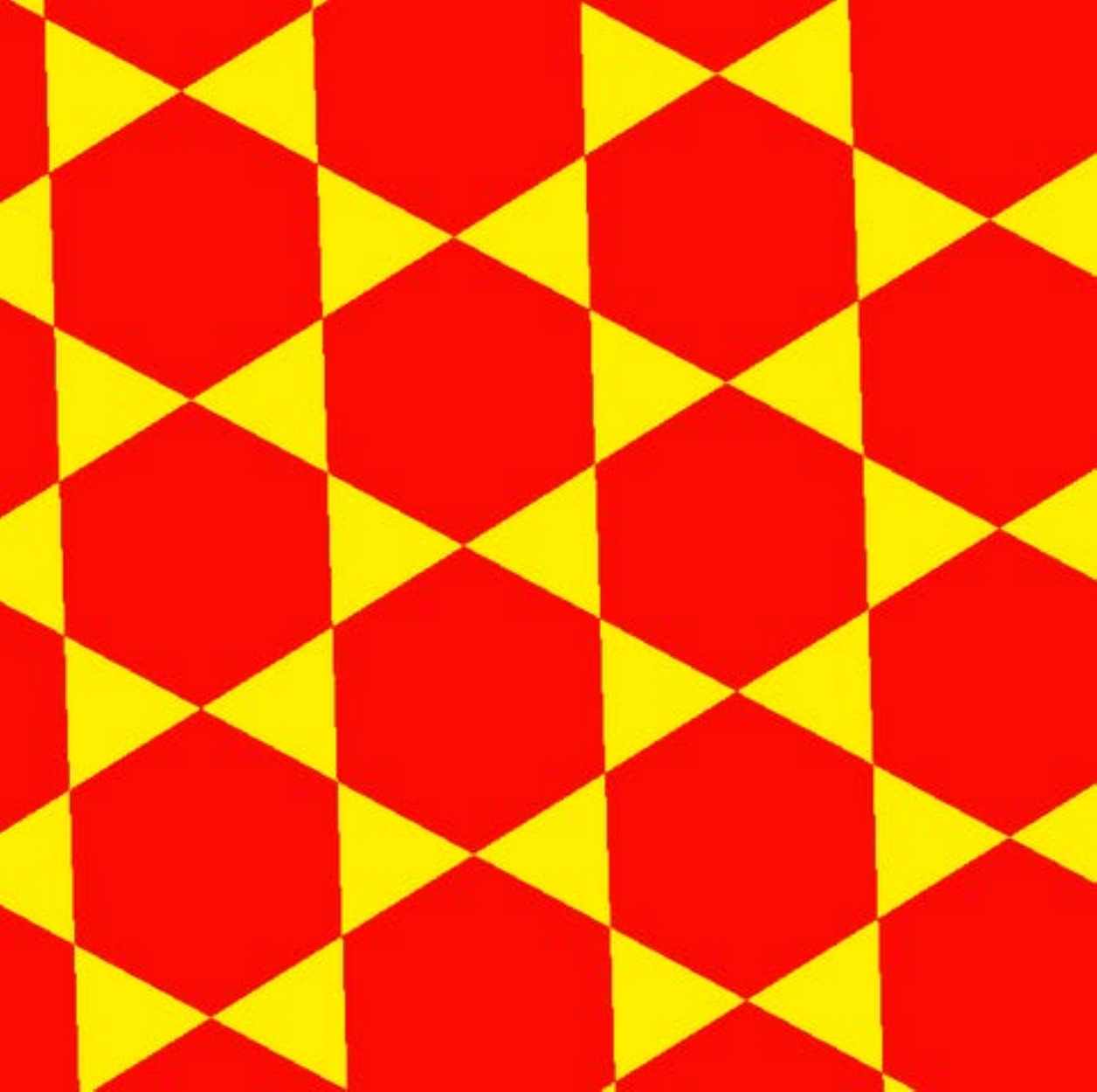} & \includegraphics[width=3cm]{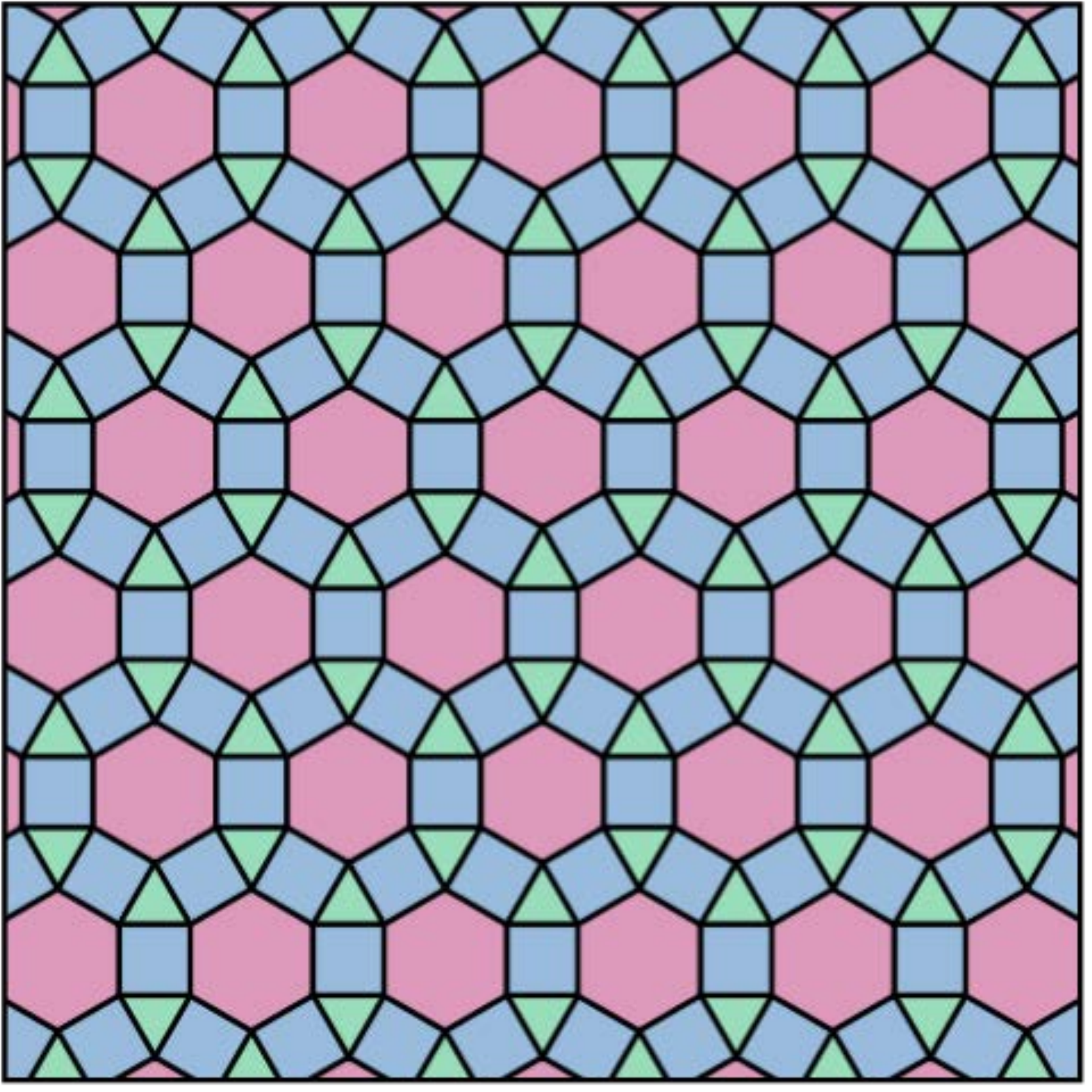} \\ \hline \hline
 \text{Name} & \text{{\small Truncated hexagonal (trH)}} & \text{Square  (S)} & \text{{\small Truncated trihexagonal (trTH)}} \\ \hline
 \text{Notation} & \left(3,12^2\right) & \left(4^4\right) & (4,6,12) \\ \hline
 \  & \includegraphics[width=3.2cm]{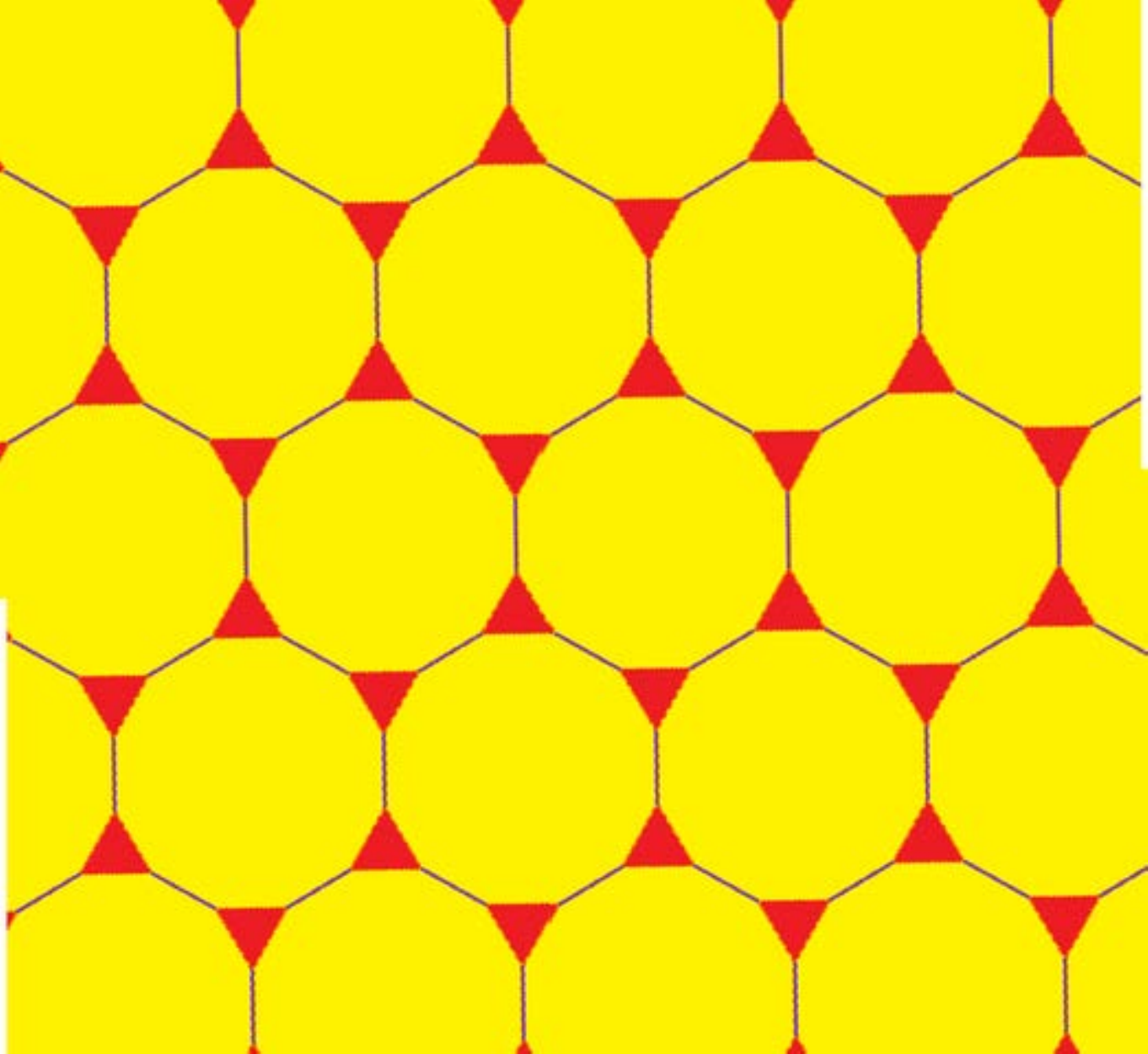} & \includegraphics[width=3cm]{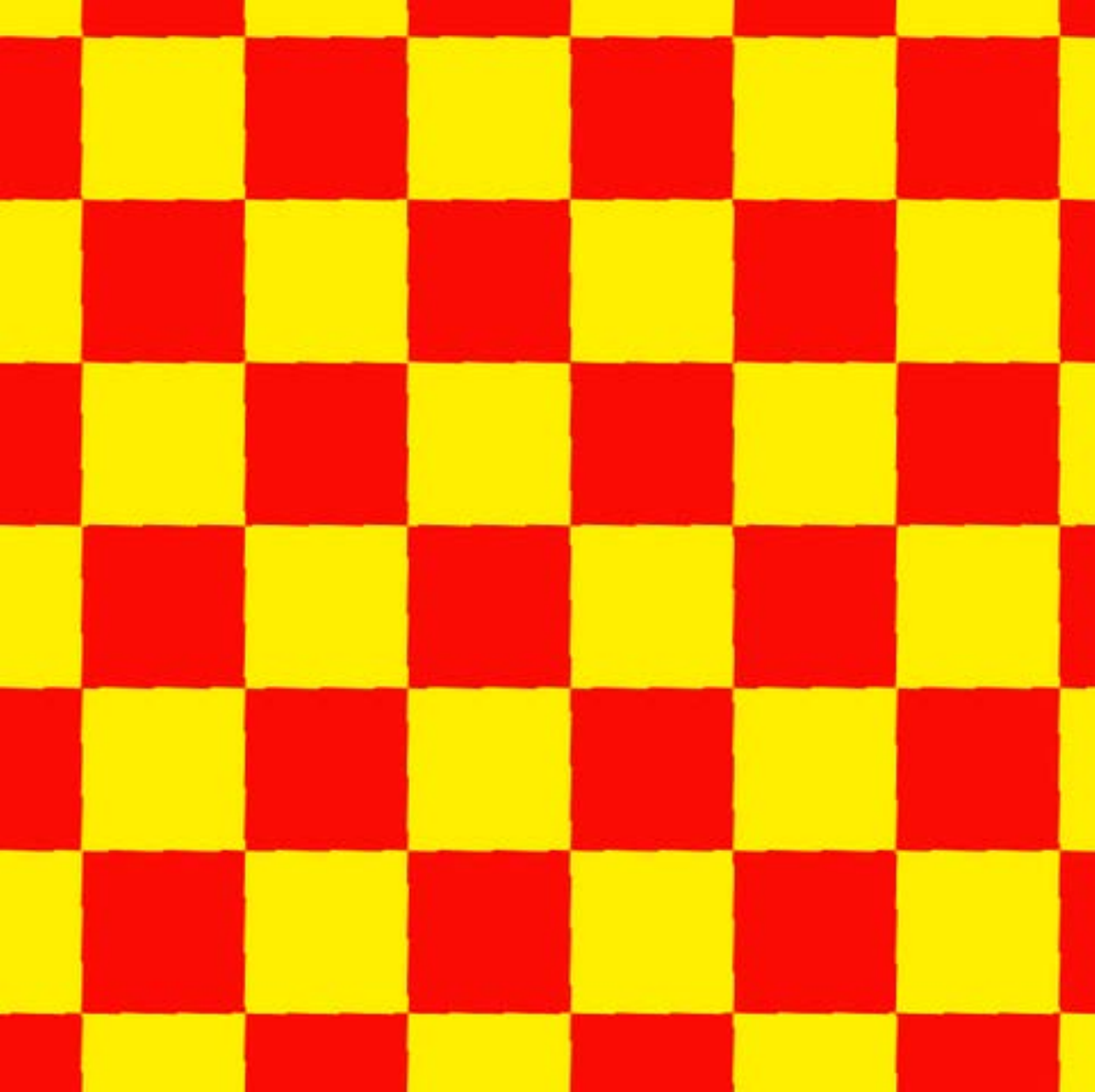} & \includegraphics[width=3cm]{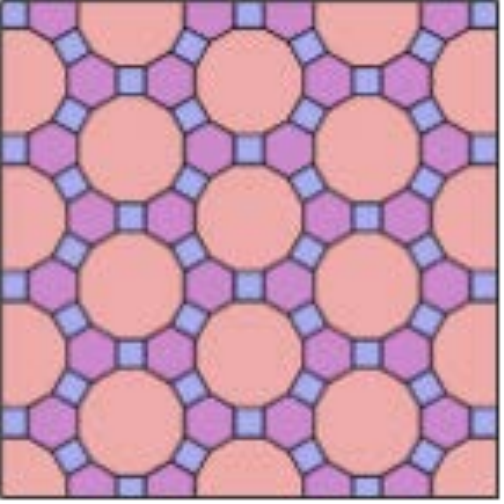} \\ \hline \hline
 \text{Name} & \text{{\small Truncated square (trS)}} & \text{Hexagonal (H)} & \ \\ \hline
 \text{Notation} & \left(4,8^2\right) & (6^3) & \ \\ \hline
 \  & \includegraphics[width=3cm]{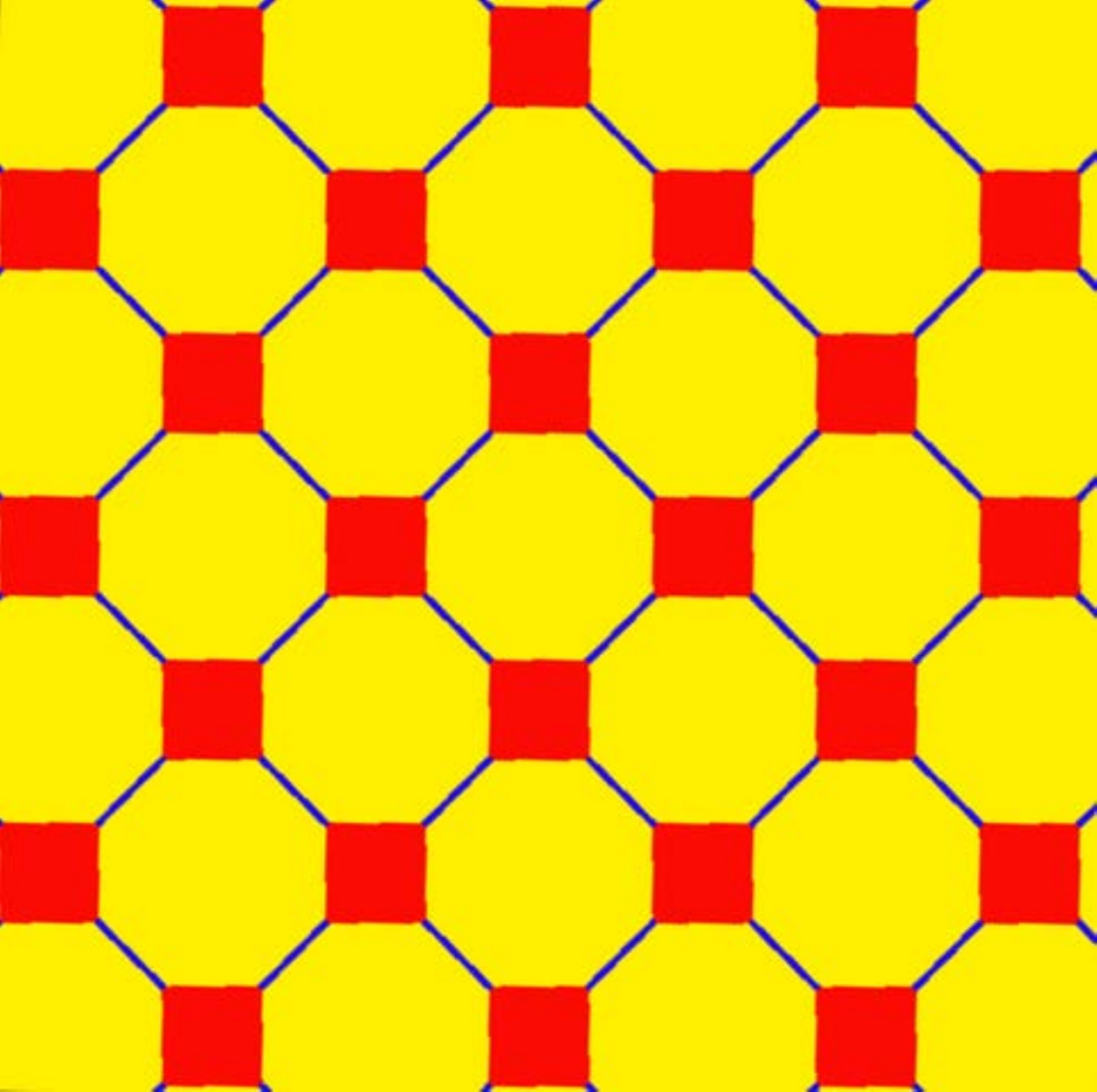} & \includegraphics[width=3cm]{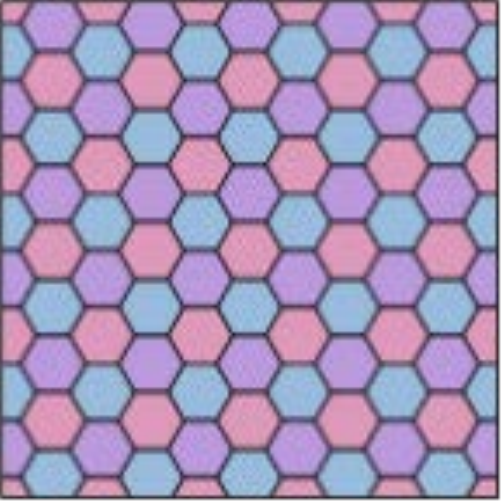} & \ \\ \hline
 \end{array}\label{Tab5.1}$
 \end{table}

  \begin{appendix}
  \section{Characteristic matrix for truncated hexagonal tiling}
 \begin{footnotesize}
 $$
 \left[
 \begin{array}{cccccccccccccccccc}
 1 & 0 & 0 & 0 & -1 & 0 & 0 & 0 & 0 & 0 & 0 & 0 & 0 & 0 & 0 & 0 & 0 & 0\\
 1 & 0 & 0 & 0 & 0 & -1 & 0 & 0 & 0 & 0 & 0 & 0 & 0 & 0 & 0 & 0 & 0 & 0\\
  0 & 0 & 0 & 0 & 0 & 0 & 0 & 0 & 0 & 1 & 0 & 0 & 0 & 1 & 1 & 0 & 0 & 0\\
   C_1 & 0 & 0 & 0 & 0 & 0 & -C_7 & 0 & 0 & S_1 & 0 & 0 & 0 & 0 & 0 & -S_7 & 0 & 0\\
   C_1 & -C_2 & 0 & 0 & 0 & 0 & 0 & 0 & 0 & S_1 & -S_2 & 0 & 0 & 0 & 0 & 0 & 0 & 0\\
   C_1' & C_2' & 0 & 0 & 0 & 0 & C'_7 & 0 & 0 & S_1' & S_2' & 0 & 0 & 0 & 0 & S_7' & 0 & 0\\
  0 & 0 & 0 & \be & 0 & 0 & -1 & 0 & 0 & 0 & 0 & 0 & 0 & 0 & 0 & 0 & 0 & 0\\
  0 & 0 & 0 & 0 & 0 & -C_6 & 1 & 0 & 0 & 0 & 0 & 0 & 0 & 0 & -S_6 & 0 & 0 & 0\\
  0 & 0 & 0 & 0 & 0 & -C_6'& 0 & 0 & 0 & 0 & 0 & 0 & \be& 0 & -S_6' & 1 & 0 & 0\\
   0 & 1 & -1 & 0 & 0 & 0 & 0 & 0 & 0 & 0 & 0 & 0 & 0 & 0 & 0 & 0 & 0 & 0\\
  0 & 1 & 0 & 0 & 0 & 0 & 0 &-1 & 0 & 0 & 0 & 0 & 0 & 0 & 0 & 0 & 0 & 0\\
  0 & 0 & 0 & 0 & 0 & 0 & 0 & 0 & 0 & 0 & 1 & 1 & 0 & 0 & 0 & 0 & 1 & 0\\
  0 & 0 & C_3 & -C_4 & 0 & 0 & 0 & 0 & 0 & 0 & 0 & S_3 & -S_4 & 0 & 0 & 0 & 0 & 0\\
  0 & 0 & C_3 & 0 & 0 & 0 & 0 & 0 & -C_9 & 0 & 0 & S_3 & 0 & 0 & 0 & 0 & 0 & -S_9\\
  0 & 0 & C_3' & C_4' & 0 & 0 & 0 & 0 & C_9' & 0 & 0 & S_3' & S_4' & 0 & 0 & 0 & 0 & S_9'\\
  0 & 0 & 0 & 0 & \al C_5 & 0 & 0 & 0 & -1 & 0 & 0 & 0 & 0 & \al S_5 & 0 & 0 & 0 & 0\\
  0 & 0 & 0 & 0 & 0 & 0 & 0 & C_8 & 0 & -1 & 0 & 0 & 0 & 0 & 0 & 0 & S_8 & 0\\
  0 & 0 & 0 & 0 & \al C_5' & 0 & 0 & C'_8 & 0 & 0 & 0 & 0 & 0 & \al S_5' & 0 & 0 & S_8'& -1
 \end{array}
 \right]
 $$
 \end{footnotesize}
\section{Some inequalities in the proof of Theorem 7.2}
This section gives the details of some tedious arguments in the proof of Theorem 7.2(c),(d) and (e)
\begin{enumerate}
	\item[(c)]
	We claim that $M_1(\theta_1,\theta_2)=\omega_1-2\omega_3-190\omega_2+191\geq 0 $.
	for any $(\theta_1,\theta_2)\in[-\pi,\pi]^2$ and $t>0$.
	Now, let $\xi=\cos\dfrac{\theta_1+\theta_2}{2}$ and $\eta=\cos\dfrac{\theta_1-\theta_2}{2}$, then we have the following:
	\begin{eqnarray*}
		\omega_2&=&\dfrac{1}{2}\xi^2+\dfrac{1}{2}\xi\eta;\\
		\omega_3&=&(2\xi^3\eta-\dfrac{3}{2}\xi\eta+\dfrac{1}{2}\eta^2); \\
		\omega_1&=&2\xi^4+2\xi^2\eta^2-3\xi^2-\eta^2+1.
	\end{eqnarray*}
	So we have
	\begin{eqnarray*}
		M_1(\theta_1,\theta_2) &=& 2(\xi^4-2\xi^3\eta+\xi^2\eta^2-46\xi\eta-49\xi^2-\eta^2+96)\\
		&=& 2 ((\xi^2-\xi\eta)^2+96-49\xi^2-46\xi\eta-\eta^2)\\
		&\geq& 0.
	\end{eqnarray*}
	\hskip0.25in
	The second part is to show that $M_2(\theta_1,\theta_2)=8\omega_1-16\omega_3+160\omega_2+37\geq 0$
	for $(\xi,\eta)\in[-1,1]^2$. Let
	\begin{eqnarray*}
		M_2(\theta_1,\theta_2)&=& 16\xi^4-32\xi^3\eta+16\xi^2\eta^2+56\xi^2+104\xi\eta-16\eta^2+45=M(\theta_1,\theta_2)\\
		&=& 16\xi^4-8\xi^2(4\xi\eta+3)+(16\xi^2\eta^2+24\xi\eta+9)+80\xi^2+80\xi\eta-16\eta^2+36\\
		&=& (4\xi^2-4\xi\eta-3)^2+20(2\xi+\eta)^2+36(1-\eta^2)\\
		&\geq& 0.
	\end{eqnarray*}
	\item[d.] Here we want to show that $M_3(\theta_1,\theta_2)=\omega_1-14\omega_3-826\omega_2+839\ge 0$
	for any $(\theta_1,\theta_2)\in[-\pi,\pi]^2$ and $t>0$. With the same change of variable as above,
	\begin{eqnarray*}
		M_3(\theta_1,\theta_2)&=& 2\xi^4-28\xi^3\eta+2\xi^2\eta^2-416\xi^2-392\xi\eta-8\eta^2+840=M(\theta_1,\theta_2)\\
		&=& 2\xi^2(\xi-\eta)^2+840-416\xi^2-392\xi\eta-8\eta^2-24\xi^3\eta\\
		&\geq& 0.
	\end{eqnarray*}
	The second part is to investigate $M_4(\theta_1,\theta_2)=8\omega_1-16\omega_3+16\omega_2+16\sqrt{3}\omega_3+48\sqrt{3}\omega_2-10\sqrt{3}+19
	\ge 0$. Consider $(\xi,\eta)\in[-1,1]^2$,
	\begin{eqnarray*}
		M_4(\theta_1,\theta_2)&=& 16\xi^4+(32\sqrt{3}-32)\xi^3\eta+16\xi^2\eta^2+(24\sqrt{3}-16)\xi^2+32\xi\eta+(8\sqrt{3}-16)\eta^2\\
		&&+27-10\sqrt{3}
	\end{eqnarray*}
	Knowing that $M_4(\th_1,\th_2)=0$ when $\d (\xi,\eta)=(\pm \frac{1}{2},\mp 1)$, we design a factorization to solve the problem.
	\begin{eqnarray*}
		M_4(\theta_1,\theta_2)&=& 16(\xi(\xi+(\sqrt{3}-1)\eta)+\frac{2\sqrt{3}-3}{4})^2+8(\sqrt{3}+1)\xi^2+40(\sqrt{3}-1)\xi\eta\\
		&& +16(2\sqrt{3}-3)\xi^2\eta^2+8(\sqrt{3}-2)\eta^2+6+2\sqrt{3}\\
		&=& 16(\xi(\xi+(\sqrt{3}-1)\eta)+\frac{2\sqrt{3}-3}{4})^2+16(2\sqrt{3}-3)(\xi^2\eta^2+\xi\eta+\frac{1}{4})\\
		&&+8(\sqrt{3}+1)(\xi^2+\xi\eta+\frac{\eta^2}{4})+6(3-\sqrt{3})(1-\eta^2)\\
		&=& 16(\xi(\xi+(\sqrt{3}-1)\eta)+\frac{2\sqrt{3}-3}{4})^2+16(2\sqrt{3}-3)(\xi\eta+\frac{1}{2})^2+8(\sqrt{3}+1)(\xi+\frac{\eta}{2})^2\\
		&& +6(3-\sqrt{3})(1-\eta^2)
	\end{eqnarray*}
	Therefore $M_4$ is nonnegative.
	\item[(e)]   Let $ t>0$. The first part is to study $M_5(\theta_1,\theta_2):=8\omega_1-272\omega_3-17648\omega_2+17912$. Observe that
	\begin{eqnarray*}
		M_5(\theta_1,\theta_2)&=& 16(\xi^4-34\xi^3\eta+\xi^2\eta^2-526\xi\eta-553\xi^2-9\eta^2+1120)=M(\theta_1,\theta_2)\\
		&=& 16 ( \xi^2(\xi-\eta)^2+1120-32\xi^2\eta-526\xi\eta-553\xi^2-9\eta^2\\
		&\geq& 0.
	\end{eqnarray*}
	So $M_5$ is nonnegative.  The second part is to show that $M_6(\theta_1,\theta_2):=16\omega_3+8\omega_1+64\omega_2-7\ge 0$. Here
	\begin{eqnarray*}
		M_6(\theta_1,\theta_2)&=& 16\xi^4+32\xi^3\eta+16\xi^2\eta^2+8\xi^2+8\xi\eta+1=M(\theta_1,\theta_2)\\
		&=& 16\xi^4+8\xi^2(4\xi\eta+1)+(4\xi\eta+1)^2\\
		&=& (4\xi^2+4\xi\eta+1)^2\\
		&\geq& 0.
	\end{eqnarray*}
	\hskip0.25in
	Finally we want to show that when $t\in(-1,1)$, for all $(\theta_1,\theta_2)\in[-\pi,\pi]^2$
	\begin{center} $g(\theta_1,\theta_2):=t^6-6t^5-9t^4+(60-48\omega_2)t^3
		+(63-48\omega_2)t^2+(144\omega_2-32\omega_3-118)t
		+8\omega_1-48\omega_3+160\omega_2-127<0, .$
	\end{center}
	\hskip0.25in
	With the same substitution as above for $\xi$ and $\eta$, we have
	\begin{align*}
	g(\theta_1,\theta_2)=& t^6-6t^5-9t^4-t^3(24\xi^2+24\xi\eta-60)-t^2 (24\xi^2+24\xi\eta-63)\\
	&-t(64\xi^3\eta-72\xi^2-120\xi\eta+16\eta^2+118)\\
	&+16\xi^4-96\xi^3\eta+ 16\xi^2\eta^2+56\xi^2+152\xi\eta-32\eta^2-119\\
	&\triangleq f(\xi,\eta)
	\end{align*}
	Hence, the problem is equivalent to showing that for a fixed $t\in(-1,1)$ then $f(\xi,\eta)<0,$ for all $(\xi,\eta)\in[-1,1]^2$.
	
	\hskip0.25in
	We shall use calculus to deal with it. To compute for the critical points of $f$ in the interior of $[-1,1]\times[-1,1]$,
	\begin{align}
	f_\xi(\xi,\eta)=&64\xi^3-192\xi^2\eta t-288\xi^2\eta+32\xi\eta^2-48\xi t^3- 48\xi t^2+144\xi t +112\xi-24\eta t^3-24\eta t^2\nonumber\\
	&+120 \eta t+152\eta\label{1}\\
	f_\eta(\xi,\eta)=&(32\xi^2-32 t-64)\eta+(-64 t-96)\xi^3+(-24 t^3-24 t^2+120 t+152)\xi
	\end{align}
	If $f_\eta(\xi,\eta)=0$ then $\eta=\dfrac{(8 t+ t^2)\xi^3+(3 t^3+3 t^2-15 t-19)\xi}{4(\xi^2- t-2)}$. Substituting this $\eta$ to (\ref{1}) we get
	\begin{equation}
	0=\dfrac{2\xi}{(\xi^2- t-2)^2}( t^2+3 t+2)(-8\xi^2+11-3 t^2)(16\xi^4+(-24 t-56)\xi^2-3 t^3-3 t^2+27 t+43)\label{2.15}
	\end{equation}
	So from (\ref{2.15}), it is clear that the critical points of $f$ are the point $\xi=0$, $\xi=\pm\sqrt{\dfrac{11-3 t^2}{8}}$, and $\xi^2=\dfrac{3 t+7\pm( t+1)\sqrt{3( t+2)}}{4}$. However, the only critical point inside $[-1,1]^2$ is at the points $\xi=0$. Moreover, $f(0,0)= t^6-6 t^5-9 t^4+60 t^3+63 t^2-118 t-119=( t+1)( t^2-2 t-7)( t^3-5 t^2-5 t+17)<0$, for all $ t\in(-1,1)$.
	
	\hskip0.25in
	Now, we focus on the boundary of $[-1,1]^2$.
	{Along} $\{1\}\times[-1,1]$:\\
	$f(1,\eta)=-( t+1)(16\eta^2+24\eta t^2-56\eta- t^5+7 t^4+2 t^3-38 t^2-2+47)$ and \\
	$f_\eta(1,\eta)=-( t+1)(24 t^2+32\eta-56)$. Thus, $f$ has no critical points.
	
	\hskip0.25in
	{Along} $\{-1\}\times[-1,1]$:\\
	$f(-1,\eta)=-( t+1)(16\eta^2-24\eta t^2+56\eta- t^5+7 t^4+2 t^3-38 t^2-2+47)$ and\\
	$f_\eta(-1,\eta)=-( t+1)(-24 t^2+32\eta+56)$. Thus, $f$ has no critical points.
	
	\hskip0.25in
	{Along} $[-1,1]\times\{1\}$:\\
	$f(\xi,1)=16\xi^4+(-64 t-96)\xi^3+(72 t-24 t^2-24 t^3+72) \xi^2+(-24 t^3-24 t^2+120 t+152)\xi+( t^6-6 t^5-9 t^4+60 t^3+63 t^2-134 t-151) $ and\\
	$f_\xi(\xi,1)=8(2\xi+1)(4\xi^2+(-12 t-20)\xi+(-3 t^3-3 t^2+15 t+19))$. It is easy to see that the critical points of $f$ are the points $\xi_1:=-\dfrac{1}{2}$ and $\xi_2:=\dfrac{3 t+5\pm( t+1)\sqrt{3( t+2)}}{2}$. One can show that $\xi_2\notin[-1,1]$, for all $  t\in(-1,1)$. Thus, $\xi_1$ is the only critical number and $f(-\dfrac{1}{2},1)=( t^2-4)( t^2-3 t-7)^2<0$, for all $ t\in(-1,1)$.
	
	\hskip0.25in
	{Along} $[-1,1]\times\{-1\}$:\\
	$f(\xi,-1)=16\xi^4+(64 t+96)\xi^3+(72 t-24 t^2-24 t^3+72) \xi^2+(24 t^3+24 t^2-120 t-152)\xi+( t^6-6 t^5-9 t^4+60 t^3+63 t^2-134 t-151) $ and\\
	$f_\xi(\xi,-1)=8(2\xi-1)(4\xi^2+(12 t+20)\xi+(-3 t^3-3 t^2+15 t+19))$. It is easy to see that the critical points of $f$ are the points $\xi_1:=\dfrac{1}{2}$ and $\xi_2:=\dfrac{-3 t-5\pm( t+1)\sqrt{3( t+2)}}{2}$. One can show that $\xi_2\notin[-1,1]$. Thus, $\xi_1$ is the only critical number and $f(\dfrac{1}{2},-1)=( t^2-4)( t^2-3 t-7)<0$.
	
	\hskip0.25in
	Also, $f(1,1)=f(-1,-1)=( t-1)^2( t+1)^3( t-7)<0$, and $f(-1,1)=f(1,-1)=( t+1)( t^2-2 t-7)( t^3-5 t^2-5 t+17)<0$. Therefore, the claim holds.
\end{enumerate}
\begin{flushright}
	$\Box$
\end{flushright}
\end{appendix}
 \end{document}